\documentclass[12pt]{article}
\usepackage{epsfig}
\usepackage{color}
\usepackage{amsfonts,amsthm,amsmath,amssymb,amscd,graphicx}

\begin{document}
\newtheorem{cor}{Corollary}
\newtheorem{lm}{Lemma}
\newtheorem{theorem}{Theorem}
\newtheorem{df}{Definition}
\newtheorem{prop}{Proposition}
\newtheorem{remark}{Remark}

\title{On the phenomenon of mixed dynamics in Pikovsky-Topaj system of coupled rotators.}
\author{Gonchenko A.S.$^1$, Gonchenko S.V.$^1$, Kazakov A.O.$^{2}$, Turaev D.V.$^{3,1,4}$ \\
$^1$ Lobachevsky State University of Nizhny Novgorod, \\
23 Prospekt Gagarina, Nizhny Novgorod, 603950, Russia \\
$^2$ National Research University Higher School of Economics, \\
25/12 Bolshaya Pecherskaya Ulitsa, 603155 Nizhny Novgorod, Russia.\\
$^3$ Department of Mathematics, Imperial College, \\
London SW7 2AZ, United Kingdom.\\
$^4$ Joseph Meyerhoff Visiting Professor, Weizmann Institute of Science \\
234 Herzl Street, Rehovot 7610001, Israel
}
\maketitle

\begin{abstract}
A one-parameter family of time-reversible systems on $\mathbb{T}^3$ is considered. It is shown that the dynamics is not conservative, namely the attractor and repeller intersect but not coincide. We explain this as the manifestation of the so-called mixed dynamics phenomenon which corresponds to a persistent intersection of the closure of stable periodic orbits and the closure of the completely unstable periodic orbits. We search for the stable and unstable periodic orbits indirectly, by finding non-conservative saddle periodic orbits and heteroclinic connections between them. In this way, we are able to claim the existence of mixed dynamics for a large range of parameter values. We investigate local and global bifurcations that can be used for the detection of mixed dynamics.
\end{abstract}

\section{Introduction}
In the present paper we study chaotic dynamics of the following system of three differential equations
\begin{equation}
\begin{array}{l}
\dot \psi_1 = 1 - 2 \varepsilon \sin \psi_1 + \varepsilon \sin \psi_2\\
\dot \psi_2 = 1 - 2 \varepsilon \sin \psi_2 + \varepsilon \sin \psi_1
+ \varepsilon \sin \psi_3\\
\dot \psi_3 = 1 - 2 \varepsilon \sin \psi_3 + \varepsilon \sin \psi_2,
\end{array}
\label{Pik1}
\end{equation}
where $\psi_i\in [0,2\pi)$, $i=1,2,3$, are angular variables, so
the phase space of (\ref{Pik1}) is a three-dimensional torus
$\mathbb{T}^3$. Note that system (\ref{Pik1}) is time-reversible: it is invariant with respect to
the reversal of time $t\to -t$ and the involution $R$:
\begin{equation}
\begin{array}{l}
\psi_1 \to \pi- \psi_3 \;\; , \;\; \psi_2 \to \pi - \psi_2 \;\; , \;\; \psi_3 \to \pi - \psi_1.
\label{involP}
\end{array}
\end{equation}

System (\ref{Pik1}) was proposed in \cite{PikTop04} by A.Pikovsky and D.Topaj as a
simple model describing dynamics of four coupled rotators.
Every individual rotator is described by the equation
$\dot\Psi_k = \omega_k$, $k=1,...,4$, where $\omega_k$ is a  constant and
$\Psi_k$ is an angular coordinate with period $2\pi$. It was assumed
that the frequencies $\omega_k$ are such that
$\omega_i-\omega_{i+1}=1, i=1,2,3$. Then, for the phase differences
$\psi_i=\Psi_i-\Psi_{i+1}$ one obtains the equations $\dot\psi_i = 1$.
The terms with $\varepsilon$ in (\ref{Pik1}) correspond to introducing a coupling
between the rotators. The particular choice of coupling in
(\ref{Pik1}) makes the system reversible, which leads to a very non-trivial dynamics.

\begin{figure}[h]
\begin{minipage}[h]{0.48\linewidth}
\center{\includegraphics[width=1\linewidth]{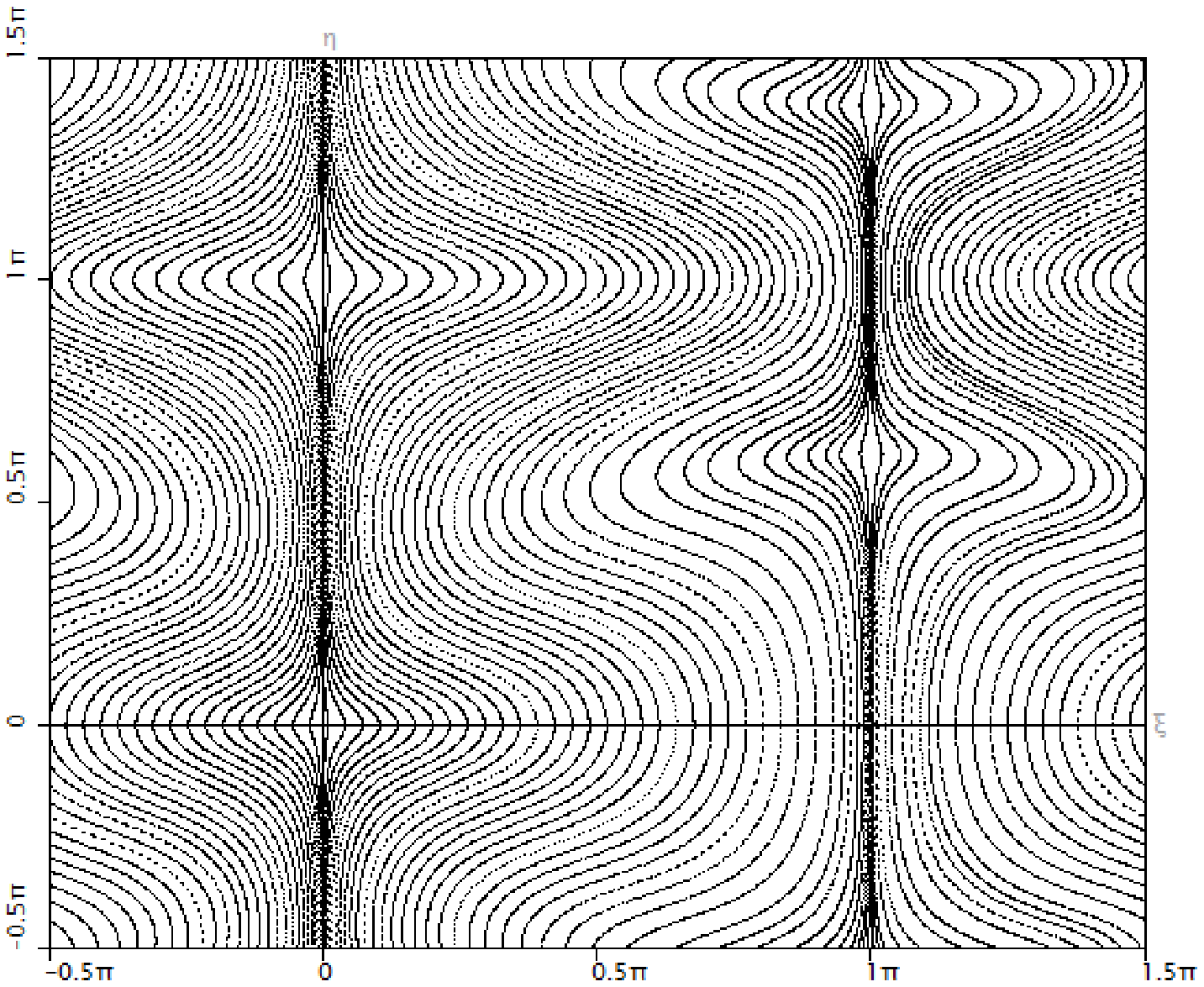} (a)}
\end{minipage}
\hfill
\begin{minipage}[h]{0.48\linewidth}
\center{\includegraphics[width=1\linewidth]{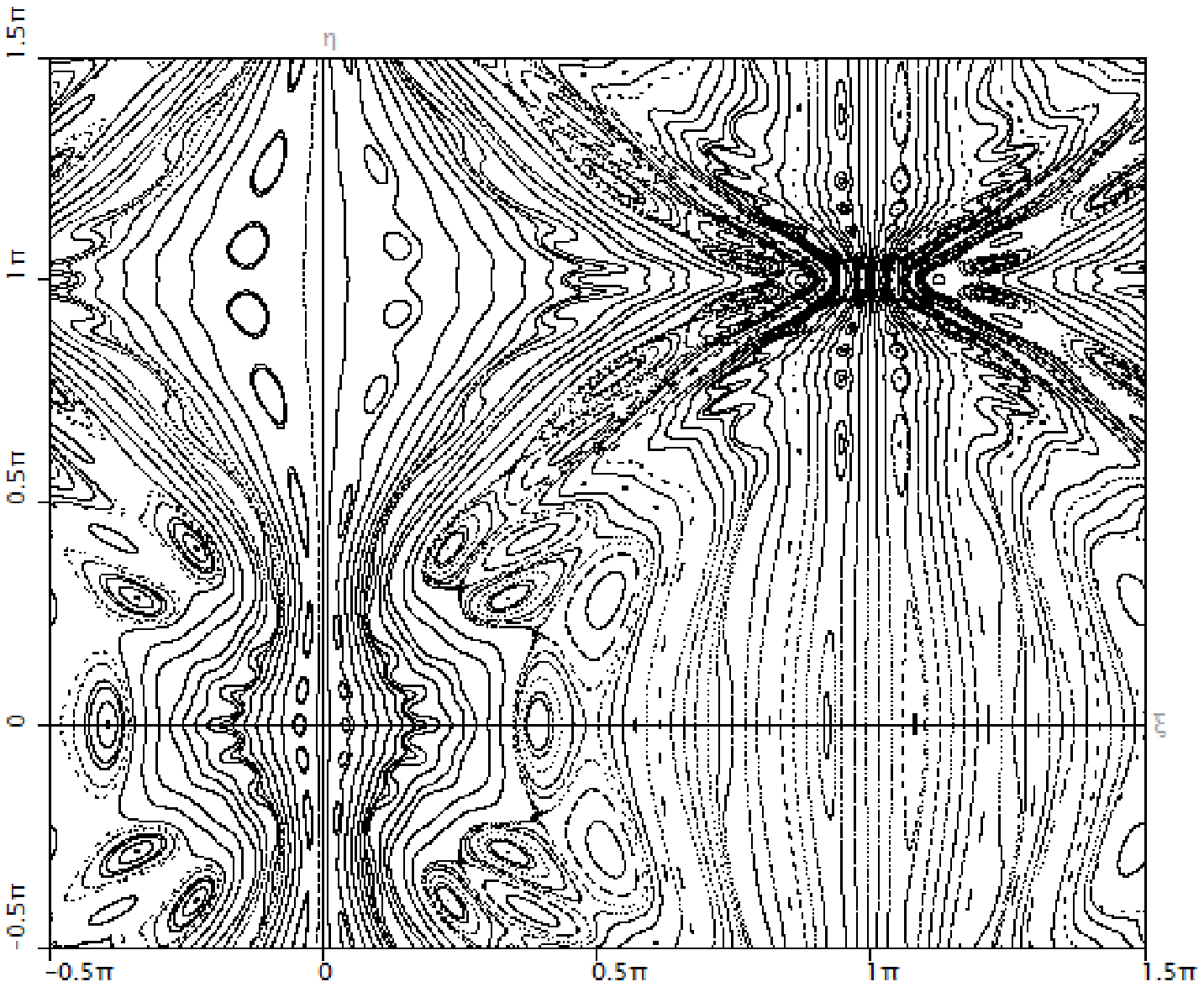} (b)}
\end{minipage}
\caption{{\footnotesize Phase portraits of the Poincar\'e map of system \eqref{Pik4} with (a) $\varepsilon = 0.1$, and (b) $\varepsilon=0.35$.
The dynamics appears conservative.}}
\label{cns1}
\end{figure}

It was noticed in \cite{PikTop04} that, at sufficiently small $\varepsilon$, the behavior of system (\ref{Pik1}) looks conservative.
In particular, for the Poincar\'e map $T_\varepsilon$ on an appropriately chosen cross-section elliptic islands are clearly observed,
see Fig~\ref{cns1}.
Moreover, the time-averaged divergence of the vector field equals to zero up to the numerical accuracy. However, with the increase
of $\varepsilon$ the apparent conservativity gets destroyed; in particular, the average divergence starts differ from zero.

\begin{figure}[h]
\begin{minipage}[h]{0.48\linewidth}
\center{\includegraphics[width=1\linewidth]{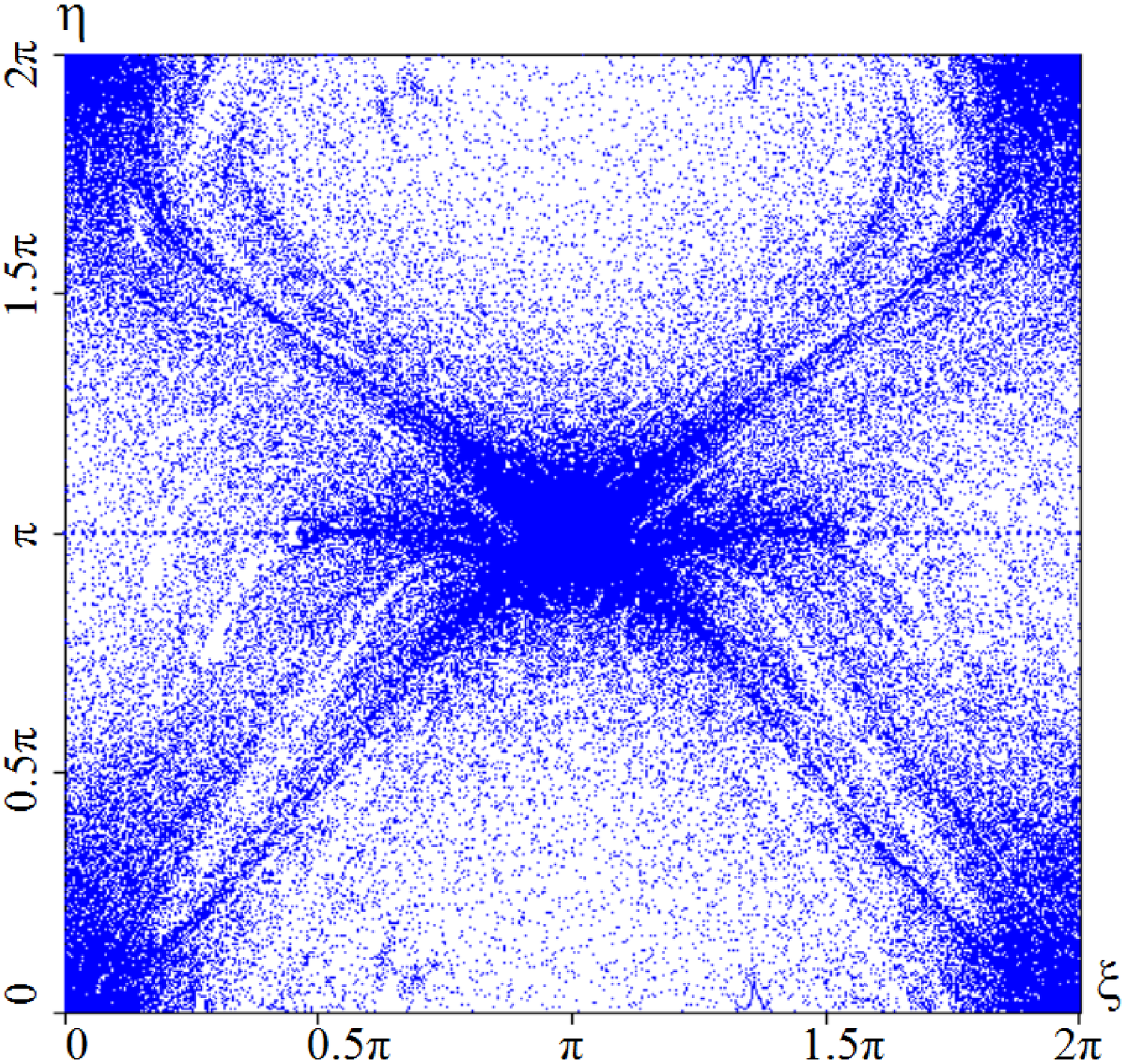} \\ (a) Reversible attractor}
\end{minipage}
\hfill
\begin{minipage}[h]{0.48\linewidth}
\center{\includegraphics[width=1\linewidth]{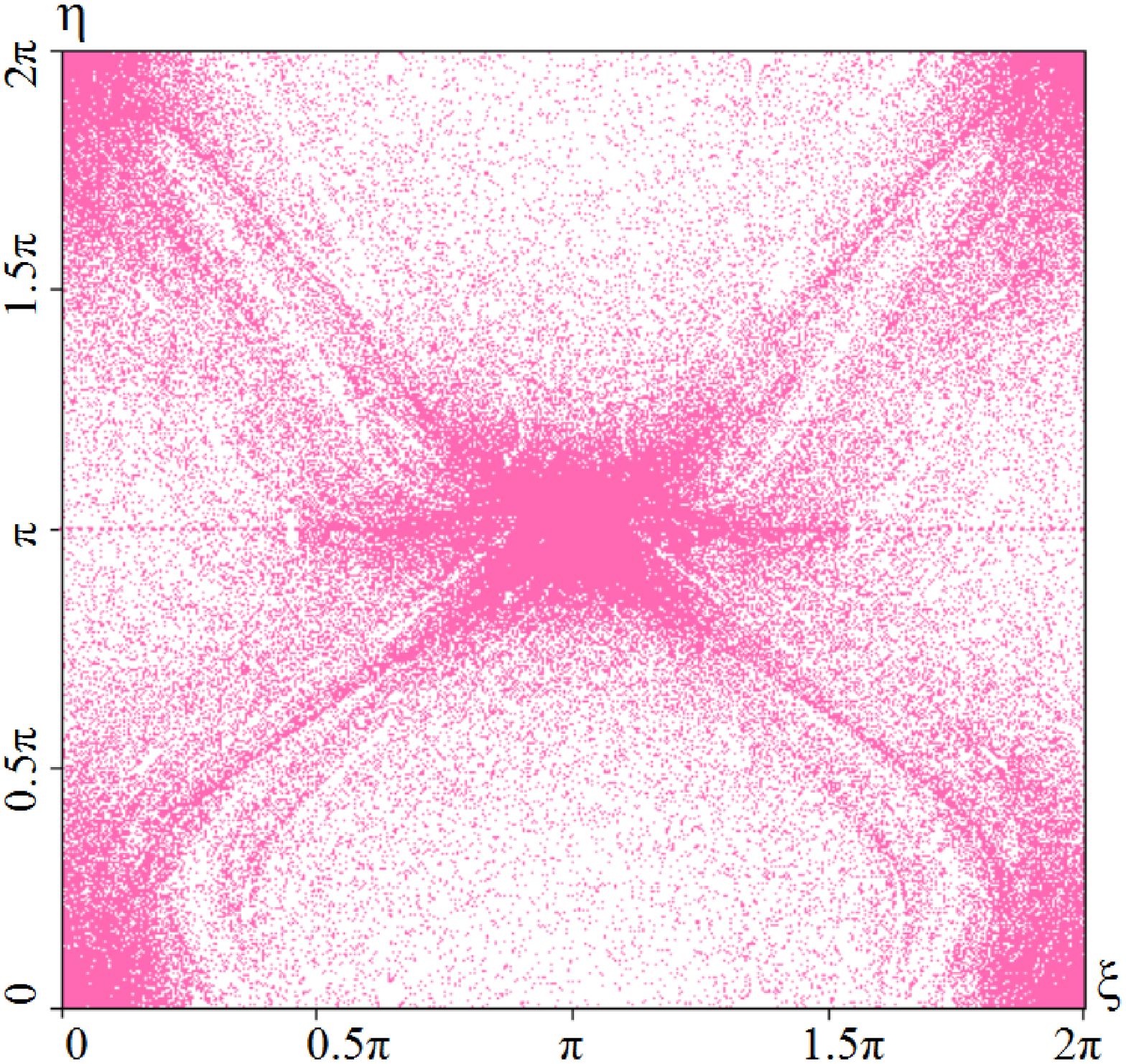} \\ (b) Reversible repeller}
\end{minipage}
\caption{{\footnotesize Phase portraits of the Poincar\'e map of system \eqref{Pik4} with $\varepsilon = 0.49$ for iterations of a uniform grid of 10000 initial conditions
on the cross-section $\Pi$; the last 100 out of 1000 iterations of each point are shown.
(a) Forward iterations; the average divergence is $\overline{\mbox{div}}\; \approx ~ -0.00122$) and (b) backward iterations $\overline{\mbox{div}}\; \approx ~ 0.00122$.
The resulting pictures are visually identical to those obtained by iterations of initial conditions uniformly distributed only
on the line $Fix(R)$ (the reversible attractor and repeller, see Section 2). Note that the numerical attractor and repeller intersect but do not coincide.}}
\label{Fig0}
\end{figure}

An interesting nonconservative effect observed in \cite{PikTop04} is the asymmetry of
the numerically obtained invariant measure for the map $T_\varepsilon$.
The invariant measure was produced as follows. On the cross-section $\psi_2=\pi/2$ take
the line $\mbox{Fix}(R): \psi_1+\psi_3 = \pi$ which
consists of the fixed points of the involution (\ref{involP}). For
a uniformly distributed set of initial points on this line, consider the iterations of
$T_\varepsilon$, and average over the iterations. The resulting sequence of point densities apparently
converges to a limit density $\mu^+$. However, beginning from a certain threshold value of $\varepsilon$, this density
is visibly non-symmetric with respect to the
involution $R$. This means that the invariant density $\mu^-$ obtained by backward iterations of the same initial points
on $Fix(R)$ is different from $\mu^+$ (these densities are related by the action of $R$),
even though the supports of these densities seem to overlap strongly.

We performed the similar computation for a uniform grid of initial conditions taken on
the cross-section $\Pi:\; \psi_1+\psi_3+2\psi_2=3\pi$ (see (\ref{transf23})-(\ref{Pik4})),
and obtained the same picture, see Fig.\ref{Fig0}. Thus , for the majority of
initial conditions the forward and backward averages differ. By Birkhoff ergodic theorem,
the forward and backward averages can be different only for a measure zero set. This means
that the numerically obtained measures $\mu^+$ and $\mu^-$ must be zero for our uniformly distributed set of initial conditions.
Therefore, the numerically produced invariant measures are, apparently, mutually singular and
not absolutely continuous with respect to the Lebesgue measure, i.e. the system is not conservative.

We interpret this phenomenon as a manifestation of
the so-called {\em reversible mixed dynamics} \cite{DGGLS13}. It corresponds to the persistent coexistence of infinitely
many periodic sinks, sources, saddles, and symmetric elliptic points. It is known
\cite{DGGLS13,GST97,LSt04,GLRT14} that for a generic reversible map the closure of the set of attracting
periodic points can intersect the closure of the set of repelling periodic points.
In this case, the numerically obtained attractor, which contains the stable periodic points and their closure, and the repeller,
which contains the unstable points and their closure, will intersect but they will not coincide.

This is very compatible with the numerical pictures for the Pikovsky-Topaj model. In this paper we show that
the reversible mixed dynamics is indeed present here. We do not search for the attracting/repelling periodic orbits directly,
as their periods are apparently very large. Instead, we establish their existence by finding non-transverse heteroclinic
cycles which include saddles of small periods (up to period 7 in our experiments). Crucially, the saddles are non-conservative, i.e., one of the saddles is
area-contracting (i.e. the Jacobian $J$ of the period map is less than 1) and the other saddle is expanding ($J>1$). It is proven in \cite{GST97} that
bifurcations of such cycles that contain both contracting and expanding saddles lead to
a simultaneous birth of infinitely many periodic attractors and repellers; see \cite{T10,T15} for generalizations to other
classes of attractors and \cite{LSt04} for the reversible case.

We find the pairs of non-conservative saddles by detecting local bifurcations of a peculiar type. We notice that the Poincare map $T_\varepsilon$ in this
model is the square of a certain orientation-reversing diffeomorphism $T_\star$ and find bifurcations which correspond to the emergence of a symmetric periodic
point of $T_\star$ with the multipliers (+1,-1). This bifurcation is described by the same normal form as the bifurcation
of periodic points with multipliers (+1,+1) with an additional symmetry (see Section \ref{sec:bnorm}).  Note that there are 4 different cases of normal forms for this
bifurcation \cite{LT12}. Two such cases have been detected in the Pikovsky-Topaj model.

The first case corresponds to the birth of an elliptic orbit and a pair of saddles, one expanding and one contracting.
These saddles are born along with heteroclinic connections, and the non-transverse intersections necessary for the proof of the mixed dynamics
appear naturally (see Fig.~\ref{Fig:Heteroclinic_large}a, \ref{Fig:FP3_heteroclinic}).

The second case corresponds to the birth of one saddle orbit, one sink, and one source. We find
this bifurcation at $\varepsilon=\varepsilon_1^*\approx 0.6042$. Because of an additional symmetry, the
Poincare map $T_\varepsilon$ has simultaneously 2 fixed points which undergo this bifurcation.
Thus, at $\varepsilon>\varepsilon_1^*$ the Poincare map $T_\varepsilon$ has 8 fixed points: 2 sinks, 2 repellers,
and 4 conservative saddles, see Fig.~\ref{Fig:FP1_heteroclinic_MD}. Most of the orbits tend to the stable fixed
points\footnote{However, at $\varepsilon < \varepsilon_1^{het} \approx 0.690$ there exist homoclinic intersections of the invariant manifolds of the saddle
fixed points. Therefore, the stable fixed points coexist with a chaotic set. Moreover, homoclinic tangencies can also exist
for such $\varepsilon$. Despite the saddle fixed points here are conservative ($J=1$), the conservativity
of the Poincare map can be violated near the orbits of tangency and, according to \cite{DGGLS13}, the reversible
mixed dynamics can exist even for some interval of $\varepsilon>\varepsilon_1^*$, although it can be hard to detect.}.
At $\varepsilon<\varepsilon_1^*$ all the fixed points disappear, and we immediately see a large chaotic
attractor (and repeller, see Fig. \ref{Fig:FP1_heteroclinic_MD2}). This phenomenon is related to the existence of homoclinic intersections of the separatrices of
the fixed point at the bifurcation moment\footnote{In a sense, this is a reversible analogue of Lukyanov-Shilnikov
bifurcation of a saddle-node with a transverse homoclinic \cite{LS78} (the so-called transition to chaos
via intermittency).}. Note that the numerically obtained attractor and repeller visibly intersect,
which means that we have a large region in the phase space corresponding to the reversible mixed dynamics.

The paper is organized as follows. In Sec.~\ref{sec:2conc} we discuss the definition of reversible attractors,
repellers, and mixed dynamics for two-dimensional reversible maps. Our approach is based on the notion
of $\varepsilon$-orbits \cite{Ruelle81,AABG85}. We also give a review of related subjects. In Sec.~\ref{sec:3} we
discuss symmetry properties of the Pikovsky-Topaj model. In Sec.~\ref{sec:bnorm} we present elements of the theory
of local symmetry-breaking bifurcations in reversible, non-orientable two-dimensional maps.
These maps naturally emerge in systems with a time-shift symmetry \cite{DKT(BCh),LR98} which is
present in the Pikovsky-Topaj model. In Sec.~\ref{sec:3a} we study
symmetry-breaking bifurcations in the model numerically. We show how these bifurcations lead to the birth of
pairs of non-conservative saddles and non-transverse heteroclinic cycles with these saddles, which creates
mixed dynamics.

\section{Mixed dynamics in two-dimensional reversible maps}  \label{sec:2conc}

A dynamical system is called {\em reversible} if it is invariant under the time reversal $t\to -t$ and a certain coordinate transformation $R$.
Obviously, $R^2$ must be a symmetry of the system; the most basic case corresponds to $R^2 = \;\mbox{Id}$, i.e. $R$ is
an involution. In the case of discrete dynamical systems (i.e., iterated maps), one says that a map $f$ is reversible, if
$f$ and  $f^{-1}$ are conjugate by $R$, i.e.,
$$f^{-1}=R\circ f \circ R.$$

A periodic orbit is called {\em symmetric} if it is invariant with respect to $R$.
Any symmetric periodic orbit possesses the following property: if it has a multiplier $\lambda$, then  $\lambda^{-1}$ is also its multiplier.
In particular, in the case of two-dimensional reversible maps, a symmetric periodic orbit can have a pair of multipliers on the unit circle,
$\lambda_{1,2}=e^{\pm i\varphi}$, where $\varphi\in (0,\pi)$, and this property
will persist for all small perturbations which do not destroy the reversibility.
Generically, the elliptic point of a sufficiently smooth reversible map is surrounded by a large set of KAM-curves \cite{Sevr}, which may make
dynamics near the elliptic point appear conservative (however, the non-conservative behavior in the resonant zones between the KAM curves
is also generic \cite{GLRT14}).

The dynamics near a non-symmetric periodic orbit of a reversible system can be arbitrary. Just note that such orbits always exist
in pairs (one orbit in the pair is mapped to the other by $R$) and the stability properties of the two orbits are opposite: the image by $R$ of
an asymptotically stable periodic orbit, a sink, is an unstable periodic orbit, a source. The image by $R$ of a non-symmetric saddle is also a saddle;
however, it is important that if at one of these saddles the Jacobian $J$ of the period map is larger than $1$, then the Jacobian $J$ is less
than $1$ at the other saddle.

\begin{figure}[ht]
\centerline{\epsfig{file=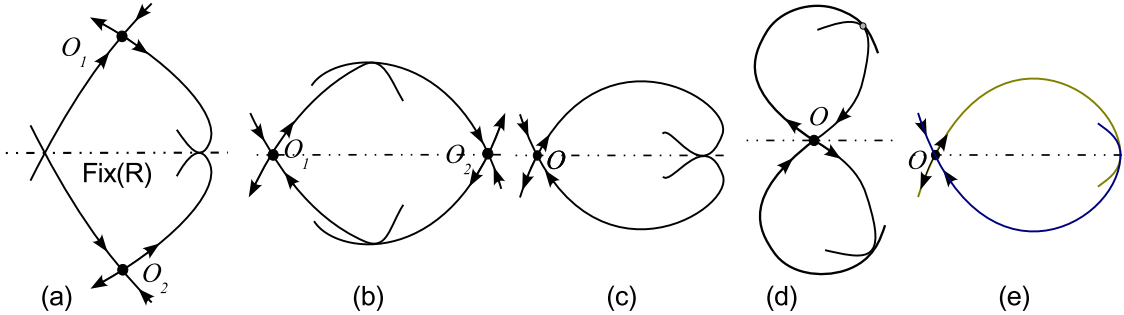, width=16cm
}} \caption{{\footnotesize Examples of two-dimensional reversible maps with symmetric homoclinic and heteroclinic tangencies.
Maps with symmetric nontransversal heteroclinic cycles are shown in (a) and (b):
here (a) $O_1 = R(O_1)$ and $J(O_1)=J(O_2)^{-1}<1$, (b) $J(O_1)=J(O_2) =1$. Case (a) is considered in \cite{LSt04}, case (b) is considered
in \cite{DGGLS13}.  Maps with symmetric homoclinic tangencies are shown in  (c)--(e): here the point $O$ is symmetric; the homoclinic orbit is symmetric in cases (c) and (e)
of a quadratic and, resp., cubic homoclinic tangency; (d) an example of a reversible map
with a symmetric pair of quadratic homoclinic tangencies to $O$. Case (d) is considered in \cite{DGGL15}.}}
\label{5typs}
\end{figure}

A remarkable phenomenon in the dynamics of reversible systems, which was called reversible mixed dynamics in \cite{DGGLS13},
is that the sinks, sources and elliptic
points can be inseparable from each other. For instance, it was shown in \cite{GLRT14} that
a generic elliptic point of a two-dimensional reversible map is a limit of a sequence of periodic sinks and a sequence of periodic sources.
The main mechanism of the emergence of the reversible mixed dynamics
is related to non-transverse heteroclinic cycles which
include a pair of non-symmetric saddles, one with $J>1$ and one with $J<1$. According to \cite{GST97,LSt04},
in a generic one-parameter unfolding of this bifurcation there exist intervals
where a typical value of the parameter corresponds to the coexistence of an infinite set
of sinks, an infinite set of sources, and an infinite set of elliptic points, and the intersection
of the closures of these sets is non-empty and contains a non-trivial hyperbolic set\footnote{This
is a generalization of the well-known ``Newhouse phenomenon'' \cite{N74,N79} to the reversible case. In terminology
of \cite{T10,T15} systems with the heteroclinic cycles of the above described type fall in the so-called ``absolute Newhouse domain''
and their dynamics is ``ultimately wild''.}.
The {\em Reversible Mixed Dynamics Conjecture} (RMD-conjecture) of \cite{DGGLS13} claims that
the same phenomenon should take place for other types of codimension-1 bifurcations of
various symmetric homoclinic and heteroclinic cycles in reversible systems. This conjecture
is proven for certain basic cases \cite{DGGLS13,LSt04,DGGL15}, see Fig. \ref{5typs}, but it
remains open in full generality, especially for the multidimensional case.

The mixed dynamics is a distinct, independent form of dynamical chaos, which should be distinguished from
the two well-known types of chaotic orbit behavior, the dissipative chaos (associated with a strange attractor)
and the conservative one (associated with a ``chaotic sea'' and elliptic islands). In our opinion,
one should expect the mixed dynamics for every non-hyperbolic two-dimensional reversible map with a
non-trivial ``symmetric dynamics''; the numerics we present in this paper supports this claim.

The main feature of the (reversible) mixed dynamics is that the attractor and the repeller of the system intersect but do not coincide.
In order to discuss this effect, we need an adequate formal definition of the attractor.
Following the ideas of \cite{Ruelle81,AABG85,Conley,TS98} we use the notion of $\varepsilon$-trajectories. Recall the definitions.

\begin{df}
Let $f: M \to M$ be a diffeomorphism of a manifold $M$ and let $\rho(x, y)$ be the
distance between points $x,y \in M$. A sequence of points $x_n\in M$ such that
$$
\rho(x_{n+1}, f(x_n)) < \varepsilon, \;\; n\in \mathbb{Z}
$$
is called an $\varepsilon$-orbit of $f$.
\label{def:epsorbit}
\end{df}

\begin{df}
We call a point $y$ attainable from a point $x$ if for any $\varepsilon >0$
there exists an $\varepsilon$-orbit which starts at $x$ and ends at $y$.
\label{def:epsdost}
\end{df}

\begin{df}
A closed invariant set $B$ is called chain-transitive, if every point of $B$
is attainable from any other point of $B$.
\label{def:epschain}
\end{df}

\begin{df}
A closed invariant set $B$ is called $\varepsilon$-stable if for every open neighborhood $U(B)$
there exists a neighborhood $V(B)$ such that, for all sufficiently small $\varepsilon>0$,
the $\varepsilon$-orbits which start in $V(B)$ never leave $U(B)$.
\label{def:epsstable}
\end{df}

\begin{df}
A chain-transitive, $\varepsilon$-stable, closed, invariant set ${\cal A}$
is called an attractor of a point $x$ if every point of $\cal A$ is attainable from $x$. A set
$\cal R$ is a repeller of a point $x$ if it is an attractor of $x$ for the inverse map $f^{-1}$.
\label{def:attrrepx}
\end{df}

It was the idea of Ruelle \cite{Ruelle81} that an attractor defined in such way would give a proper
picture of behavior of a system subject to a bounded noise. In particular, this notion is
convenient for the analysis of numerical experiments.

Note that a point may have several attractors by this definition. For any given $\delta>0$,
for all sufficiently small $\varepsilon>0$ a typical $\varepsilon$-orbit of a point $x$ will, eventually,
enter the $\delta$-neighborhood of any of the attractors of $x$ and will never leave it. Moreover,
it will then visit any neighborhood of any point of this attractor infinitely many times.
In the simplest case, if a point $x$ belongs to the domain of attraction of some periodic sink $p_s$,
then $p_s$ is the only attractor of $x$. In the reversible situation, the source $R(p_s)$ will be
the only repeller of every point in its neighborhood.

Given a set $C$, the union of all attractors of all points of $C$ will be called the attractor ${\cal A}_C$ of this set
and the union of all repellers of all points of $C$ will be called the repeller ${\cal R}_C$ of this set.

\begin{df}
Let $f$ be an $R$-reversible two-dimensional diffeomorphism, and $\dim Fix(R) = 1$.  The sets ${\cal A} = {\cal A}_{Fix{R}}$
and ${\cal R} = {\cal R}_{Fix{R}}$ are called a reversible attractor and a reversible repeller of the map $f$.
\label{def:revatrrep}
\end{df}

Evidently, $R({\cal A})= {\cal R}$. The case when ${\cal A}$ and ${\cal R}$ do not intersect
is easy to imagine: just let ${\cal A}$ lie on one side of $Fix(R)$,
then ${\cal R}$ will lie on the other side. In this case we will have a usual dissipative dynamics.
If $f$ is an area-preserving map of a compact manifold $M$, then the only chain-transitive set is the whole phase space $M$, so
${\cal A} = {\cal R}=M$. The most interesting case, where ${\cal A}\cap {\cal R}\neq \emptyset$ and ${\cal A} \neq {\cal R}$,
is different from both the dissipative and conservative cases. Such picture was first observed in \cite{Politi} and explained as the ``coexistence of conservative chaos
with the dissipative behavior'', see also \cite{LR98}. We, however, do not interpret the intersection of the attractor and repeller as a conservative
phenomenon, because the existence of non-transverse homoclinic and heteroclinic intersections
implies the birth of sinks and sources, i.e. mixed dynamics, as it is explained above and demonstrated below.
See also more discussion in \cite{G16}.

\section{Symmetries in the model.} \label{sec:3}
Let us consider system \eqref{Pik1}.
By means of the coordinate change
\begin{equation}
\begin{array}{l}
\xi=\frac{\displaystyle \psi_1 - \psi_3}{\displaystyle 2}, \;\;  \eta=\frac{\displaystyle \psi_1 +
\psi_3 - \pi}{\displaystyle 2}, \;\;
\rho=\frac{\displaystyle \psi_1 + \psi_3 - \pi}{\displaystyle 2}+\psi_2-\pi,\\
\end{array}
\label{transf23}
\end{equation}
the system is brought  to the following form
\begin{equation}
\begin{array}{l}
\dot \xi = 2 \varepsilon \sin \xi \sin \eta,\\
\dot \eta = 1 - \varepsilon \cos(\rho-\eta)  - 2 \varepsilon \cos \xi \cos \eta,\\
\dot \rho = 2 + \varepsilon \cos (\rho -\eta).
\end{array}
\label{Pik3}
\end{equation}
After the time change  $d t_{new}=(2+\varepsilon \cos(\rho-\eta))dt$
system (\ref{Pik3}) recasts as
$$
\begin{array}{l}
\dot \xi = \frac{\displaystyle 2 \varepsilon \sin \xi \sin \eta}
{\displaystyle2 + \varepsilon \cos (\rho -\eta)},\\
\dot \eta = \frac{\displaystyle 1 - \varepsilon \cos(\rho-\eta)  - 2 \varepsilon \cos \xi \cos
\eta}
{\displaystyle 2 + \varepsilon \cos (\rho -\eta)}, \\
\dot \rho = 1,
\end{array}
$$
i.e., a non-autonomous time-periodic system
\begin{equation}
\begin{array}{l}
\dot \xi = \frac{\displaystyle 2 \varepsilon \sin \xi \sin \eta}
{\displaystyle2 + \varepsilon \cos (t -\eta)},\\
\dot \eta = \frac{\displaystyle 1 - \varepsilon \cos(t-\eta)  - 2 \varepsilon \cos \xi \cos
\eta}
{\displaystyle 2 + \varepsilon \cos (t -\eta)}. \\
\end{array}
\label{Pik4}
\end{equation}
Note that system \eqref{Pik4} is well-defined for all $\varepsilon < 2$.

System \eqref{Pik4} as well its time-shift maps possess symmetries of various types.
First of all, we note the following simple facts.

\begin{itemize}
\item[(i)]
System \eqref{Pik4} is reversible, i.e., invariant with respect to the involution
\begin{equation}
R: \;\;
\begin{array}{l}
\xi \to \xi,  \;\; \eta \to -\eta
\end{array}
\label{invol}
\end{equation}
and the time reversal $t \to -t$.
\item[(ii)]
System \eqref{Pik4} is invariant with respect to the coordinate change
\begin{equation}
\sigma: \;\;
\begin{array}{l}
\xi \to \pi - \xi,  \;\; \eta \to \pi + \eta
\end{array}
\label{invol_s}
\end{equation}
and the time shift $t \to t + \pi$.
\end{itemize}

Property (ii) is called the {\em time-shift symmetry} \cite{LR98}. It is often met in problems where autonomous systems with periodic perturbations
are considered, see e.g. \cite{DKT(BCh),LR98}. Denote as $T_{a\to b}$ the time-shift map along orbits of system \eqref{Pik4}
from $t=a$ to $t=b$. It is easy to see that
\begin{equation}
T_{0 \to 2\pi} = (\sigma T_{0 \to \pi})^2.
\label{eq:T2S}
\end{equation}
Indeed, property (ii) implies that $T_{\pi \to 2\pi} = \sigma^{-1} T_{0 \to \pi} \sigma$ and, since $\sigma = \sigma^{-1}$, we have
$$
T_{0 \to 2\pi} = T_{\pi \to 2\pi} T_{0 \to \pi} =  \sigma^{-1} T_{0 \to \pi} \sigma T_{0 \to \pi} = \sigma T_{0 \to \pi} \sigma T_{0 \to \pi} = (\sigma T_{0 \to \pi})^2.
$$

In what follows, we will use the notation $T$ for the Poincar\'e
map $T_{0 \to 2\pi}$ and $T_\star$ for its square root $\sigma T_{0 \to \pi}$.
Note that the map $\sigma$ is orientation reversing, while the time-shift map $T_{0 \to \pi}$ preserves the orientation. Thus, we
obtain an interesting fact that the Poincar\'e map $T$ for system \eqref{Pik4} is the second iteration of an orientation reversing map $T_\star$.

Importantly, the map $T_\star$ is reversible with respect to the involution $R$. To show this, we need to check
the following relations:
\begin{equation}
R(\sigma T_{0 \to \pi}) = (\sigma T_{0 \to \pi})^{-1}R = (T_{0\to\pi})^{-1} \sigma^{-1} R = (T_{0\to\pi})^{-1} \sigma R.
\label{invol_sr}
\end{equation}
By \eqref{invol} and \eqref{invol_s}, we have
$$
R\sigma = \{\bar\xi = \pi - \xi, \bar\eta = \pi -\eta \},
$$
since $\xi$ and $\eta$ are $2\pi$-periodic coordinates. Analogously,
$$\sigma R  = \{\bar\xi = \pi - \xi, \bar\eta = \pi - \eta \}.
$$
Thus, relation \eqref{invol_sr} will be proven, if we show that the map  $T_{0 \to \pi}$ is reversible with respect to the involution $R\sigma=\sigma R$.
Now note that this just follows from the two facts: system \eqref{Pik4} is invariant with respect to the transformation
$\xi \to \pi - \xi, \eta \to \pi - \eta, t \to \pi - t$, and $(T_{0\to\pi})^{-1} = T_{\pi\to 0}$.

We also remark that system (\ref{Pik4}) possesses an additional symmetry $S: \xi\mapsto -\xi$  or, since (\ref{Pik4}) is given on a torus,
$\xi\mapsto 2\pi -\xi$. This property is evidently inherited by the Poincare map $T = T_{0\to 2\pi}$, so the phase portrait of $T$
is symmetric relative to the reflection with respect to the lines $\xi=0$ and $\xi =\pi$. Note also that these symmetry lines are invariant with respect to the
map $T$ and divide the torus into two invariant annuli, $0\leq \xi\leq \pi$ and $\pi \leq \xi\leq 2\pi$. The annuli are also invariant with respect to $T_\star$,
however the map $T_\star$ maps the line $\xi=0$ to $\xi=\pi$, and vice versa.

\section{Local conservativity breaking bifurcations in orientation reversing maps.}
\label{sec:bnorm}

The periodic orbits of the Pikovsky-Topaj model correspond to periodic orbits of an orientation-reversing map $T_\star$, so we need
to recall the bifurcation theory for such maps, see \cite{LT12} for a more detailed account. For an orbit of an odd period $q$, the period map
$T_\star^q$ is orientation reversing, i.e. bifurcations of such orbits are described by the theory of bifurcations of fixed points of orientation-reversing maps.
Bifurcations of orbits of even periods for $T_\star$ are described by the theory of bifurcations of fixed points
of orientation-preserving maps in general. However, in our model certain periodic orbits (e.g. the fixed points of $T=T_\star^2$ and some fixed points of $T^3=T_\star^6$)
are $S$-symmetric, i.e., they lie on the invariant lines $\xi=0$ and $\xi=\pi$. For the bifurcation of the birth of the $S$-symmetric points
the corresponding normal form is the same as in the orientation-reversing case, i.e. the theory below is also applicable to them too.

Let us consider a one-parameter family $f_\nu$ of two-dimensional maps. Assume that
the maps are reversible, i.e. there exists an involution $R$ such that $f_\nu^{-1}=R\circ f_\nu\circ R$.
Assume that the maps are orientation-reversing, i.e. $det f'_\nu <0$.

Let at $\nu=0$ the map have an $R$-symmetric fixed point $O$, i.e. $f_0 O=O$ and $RO=O$.
We assume that the set $Fix (R)$ of fixed points of the involution $R$ is one-dimensional. Then, by Bochner theorem \cite{Bochner},
one can choose coordinates near $O$ such that the involution $R$ is locally given by
\begin{equation}
R: \left\{\begin{array}{l} x \to x \\ y \to -y
\end{array}
\right. \label{involR}
\end{equation}
for all sufficiently small values of the parameter $\nu$. Thus, $Fix (R)$ is the line $y=0$.

Let $\lambda_{1,2}$ be the multipliers of $O$. By the reversibility, $\lambda_{1,2}^{-1}$ must also be the multipliers.
Since $det f'_0(0)=\lambda_1\lambda_2<0$, the only option is that
$$
\lambda_1 = -1, \;\; \lambda_2 = +1.
$$
The corresponding eigenvectors must be either orthogonal or parallel to the line $Fix (R)$. We assume that
the eigenvector corresponding to the multiplier $+1$ is orthogonal to $Fix (R)$.\footnote{The second variant (the eigenvector corresponding to $-1$ is
orthogonal to  $Fix (R)$) would correspond to a very degenerate case: here the map $R\circ f_0$ would be an involution with the linear part equal to identity;
by Bochner theorem an involution is conjugate to its linear part; the only map conjugate to identity is the identity map itself, so $R\circ f_0 =id$ in this case,
i.e. $f_0 \equiv R$ - this is definitely not the case of our model.}.

Put the fixed point $O$ at the origin of the $(x,y)$-plane. The Taylor expansion for the map $f_0$ at $O$ is as follows:
\begin{equation}
f_0 \; : \; \left\{
\begin{array}{l}
\bar x = -x + A_{20}x^2 + A_{11}xy + A_{02} y^2 + o(x^2+y^2),\\
\bar y = y + B_{20} x^2 + B_{11}xy + B_{02} y^2 + o(x^2+y^2).\\
\end{array}
\right. \label{f0}
\end{equation}
Note that the reversibility does not impose any restriction on the coefficients of the quadratic terms, i.e. the inverse map (after the change $x\to x, y\to -y$) will
have the same form (\ref{f0}) up to higher order terms. The coordinate transformation
$$x_{new}=x+ax^2+by^2, \qquad y_{new}=y+cxy$$
keeps the map $R$-reversible and brings it to the form
\begin{equation}
f_0 \; : \;\left\{
\begin{array}{l}
\bar x = -x + A_{11} xy + \dots,\\
\bar y = y + B_{20} x^2 + B_{02} y^2 + \dots,\\
\end{array}
\right.
\end{equation}
if we choose $a=A_{20}/2, b=A_{02}/2, c=-B_{11}/2$. We further assume that the following genericity condition holds:
$$A_{11} \neq 0, B_{20} \neq 0, B_{02} \neq 0.$$

By the scaling $x \to x \frac{1}{\sqrt{|B_{20}B_{02}|}},  \;\; y \to y \frac{1}{|B_{02}|}$,
we bring the map to the form
$$\left\{
\begin{array}{l}
\bar x = -x - \alpha xy + \dots\\
\bar y = y \pm x^2 \pm  y^2 + \dots\\
\end{array}
\right.$$
where $\alpha =- A_{11}/|B_{02}|$. Note that every combination of the signs ``$+$'' ``$-$'' are possible, but the cases $(-,-)$ and
$(-,+)$ can be reduced to $(+,+)$ and $(+,-)$, respectively, by considering the inverse map $f_0^{-1}$ instead $f_0$.
Thus, we finally obtain the following normal form:
\begin{equation}
f_0 \; : \;\left\{
\begin{array}{l}
\bar x = -x - \alpha xy + \dots\\
\bar y = y + x^2 \pm  y^2 + \dots\\
\end{array}
\right. \label{nf_f0}
\end{equation}

In order to study the bifurcations of the zero fixed point, we need to consider a generic one-parameter
unfolding within the class of reversible maps. We refer the reader to \cite{LT12} where it is shown that
the generic one-parameter family of perturbations of $f_0$ which keeps the map reversible can be written in the form
\begin{equation}
f_\nu \left\{
\begin{array}{l}
\bar x = (-1 - \alpha \frac{\nu}{2})x - \alpha xy + \dots,\\
\bar y = \nu + (1 \pm \nu)y + x^2 \pm  y^2 + \dots .\\
\end{array}
\right. \label{nfth}
\end{equation}
where the dots stand for cubic and higher order terms in $(x,y,\nu)$.
Up to the terms of the third order and higher, map (\ref{nfth}) coincides with the
composition of two maps: the symmetry $S: x \to -x, \; y \to y$ and the time-1 map of the
two-dimensional, reversible, autonomous flow
\begin{equation}
\begin{array}{l}
\dot x = \alpha xy,\\
\dot y = \nu + x^2\pm y^2.
\end{array}
\label{potok}
\end{equation}

\begin{figure} [ht]
\centerline{\epsfig{file=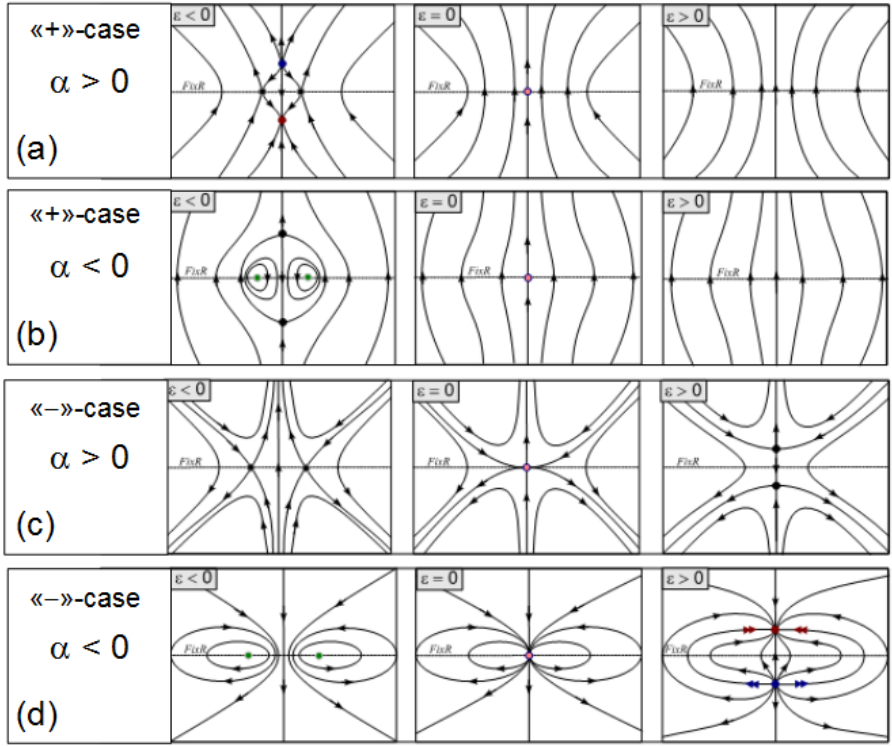,height=10cm}} \label{bifneor}
\caption{{\footnotesize Bifurcations of equilibria in the flow normal form (\ref{potok}).}}
\label{4typlocbif}
\end{figure}

Equilibrium states of system (\ref{potok}) belong either to the line $x=0$ ($R$-asymmetric equilibria) or to the line $y=0$ ($R$-symmetric
equilibria). In the ``$+$'' case, there are no equilibria at $\nu<0$ and 4 equilibria at $\nu>0$; in the ``$-$'' case we have two equilibria
at $\nu<0$ ($R$-symmetric) and two $R$-asymmetric equilibria at $\nu>0$. So, we will speak of the $0\to 4$ and $2\to 2$ bifurcations;
the corresponding bifurcation diagrams are presented in Fig.~\ref{4typlocbif}.

Note that the flow of (\ref{potok}) commutes with the map $S$. Hence, $S$-symmetric equilibria correspond to fixed points of the map
(\ref{nfth}), while $S$-symmetric pairs of $S$-asymmetric equilibria correspond to orbits of period 2. As $S$-symmetric equilibria happened to be $R$-asymmetric here,
and vice versa, we obtain the following description of the bifurcations in the orientation-reversing maps $f_\nu$.
\begin{description}
\item $\bullet$ In the ``$+$'' case, the degenerate fixed point that exists at $\nu=0$
splits into 4 periodic points at $\nu>0$: an $R$-symmetric elliptic orbit of period 2 and two $R$-asymmetric saddle fixed points
if $\alpha < 0$, and an $R$-symmetric saddle orbit of period 2 and two $R$-asymmetric fixed points, a sink and a source, if $\alpha>0$.

\item $\bullet$ In the ``$-$'' case, at $\nu > 0$ there are two $R$-asymmetric fixed points, saddles if $\alpha>0$ or a sink-source pair if $\alpha<0$.
At $\nu=0$ they merge into the single degenerate fixed point $O$, which becomes, at $\nu<0$, an $R$-symmetric orbit of period 2, saddle if $\alpha<0$ and
elliptic if $\alpha>0$.
\end{description}

As we mentioned, the same normal form (\ref{potok}) describes bifurcations of $S$-symmetric fixed points of an orientation-preserving reversible map, which
have a pair of multipliers equal to $1$. Namely, the normal form for such map coincides with the time-1 flow of (\ref{potok}) up to the terms
of the third order and higher. Thus, we have the behavior similar to the orientation-reversing case, with the difference that the $R$-symmetric period 2 orbits
become pairs of $R$-symmetric fixed points ($S$-symmetric to each other) \cite{LT12}.

Recall that in the Pikovsky-Topaj model periodic points for the Poincare map $T$ are periodic
points for the orientation-reversing map $T_\star$. If there is an orbit of period $q$ for $T$, then either it corresponds to an orbit of the same period for $T_\star$
(this can happen only if $q$ is odd), or there is one more $q$-periodic orbit of $T$ (the image of the first one by the map $T_\star$) such that together
these two orbits form an orbit of period $2q$ for the map $T_\star$. In the first case, bifurcations are described by the above theory for orientation-reversing maps
(just note that a period-2 point for $T_\star^q$ is a fixed point for $T^q$), while in the second case bifurcations are described by the theory for orientation-preserving
maps - a particular case of this theory corresponding to $S$-symmetric periodic orbits is again described by Fig.~\ref{4typlocbif}. We have encountered both situations
(orientation-reversing and $S$-symmetric orientation-preserving) in the model. In all cases we have found only bifurcations of $0\to 4$ type, with $\alpha>0$
for bifurcations of the fixed points of $T$ and $\alpha<0$ for the bifurcations of points of period 3 and 7 for the map $T$.

Note that the $0\to 4$ bifurcation with $\alpha<0$ produces a pair of $R$-asymmetric saddles, whose stable and unstable manifolds coincide for system (\ref{potok}).
As the time-1 shift by the flow of (\ref{potok}) is only an approximation to the map $f_\nu$ (up to multiplication to $S$), these separatrix connections will
split for a generic family $f_\nu$ and transverse heteroclinic orbits will form. In the $S$-symmetric case, however, one of the separatrices of each saddle must coincide
with the symmetry line, so these separatrices cannot split (while the other pair of separatrices splits). See Fig.~\ref{Fig:Heteroclinic_large}a for examples. Note that the splitting of heteroclinic
connections leads also to a creation of non-degenerate homoclinic tangencies at certain parameter values; this phenomenon plays a central role in our analysis
of the emergence of mixed dynamics.

\section{Local and global symmetry-breaking bifurcations. Numerical experiments.} \label{sec:3a}

In this section we report the results of numerical experiments that demonstrate
both local and global bifurcations of breaking reversible symmetry
in the dynamics of Poincar\'e map $T_\varepsilon$ for the Pikovsky-Topaj model.
We detect the moment of the birth of chaotic mixed dynamics which exists alongside with a simple periodic attractor and repeller,
and trace the existence of mixed dynamics from sufficiently large to sufficiently small values of $\varepsilon$
when the dynamics appears practically indistinguishable from conservative.

Note that the points of intersection of the line $Fix R$ with $T^k(FixR)$ correspond to $R$-symmetric periodic orbits of $T$.
Indeed, if $x \in Fix R$ and $T^k(x) \in Fix R$, then $R(x) = x$ and $R(T^k(x)) = T^k(x)$, then the
identity $T^{-k}R(x) = RT^k(x)$ gives that $T^{-k}(x) = T^k(x)$, i.e $T^{2k}(x) = x$.
A tangency between $Fix R$ with $T^k(Fix R)$ corresponds to a change in the number of periodic points as the parameters
of the system change, so we use the tangencies in order to detect bifurcation moments of the birth of new $R$-symmetric periodic points.

\begin{figure}[h]
\begin{minipage}[h]{0.48\linewidth}
\center{\includegraphics[width=1\linewidth]{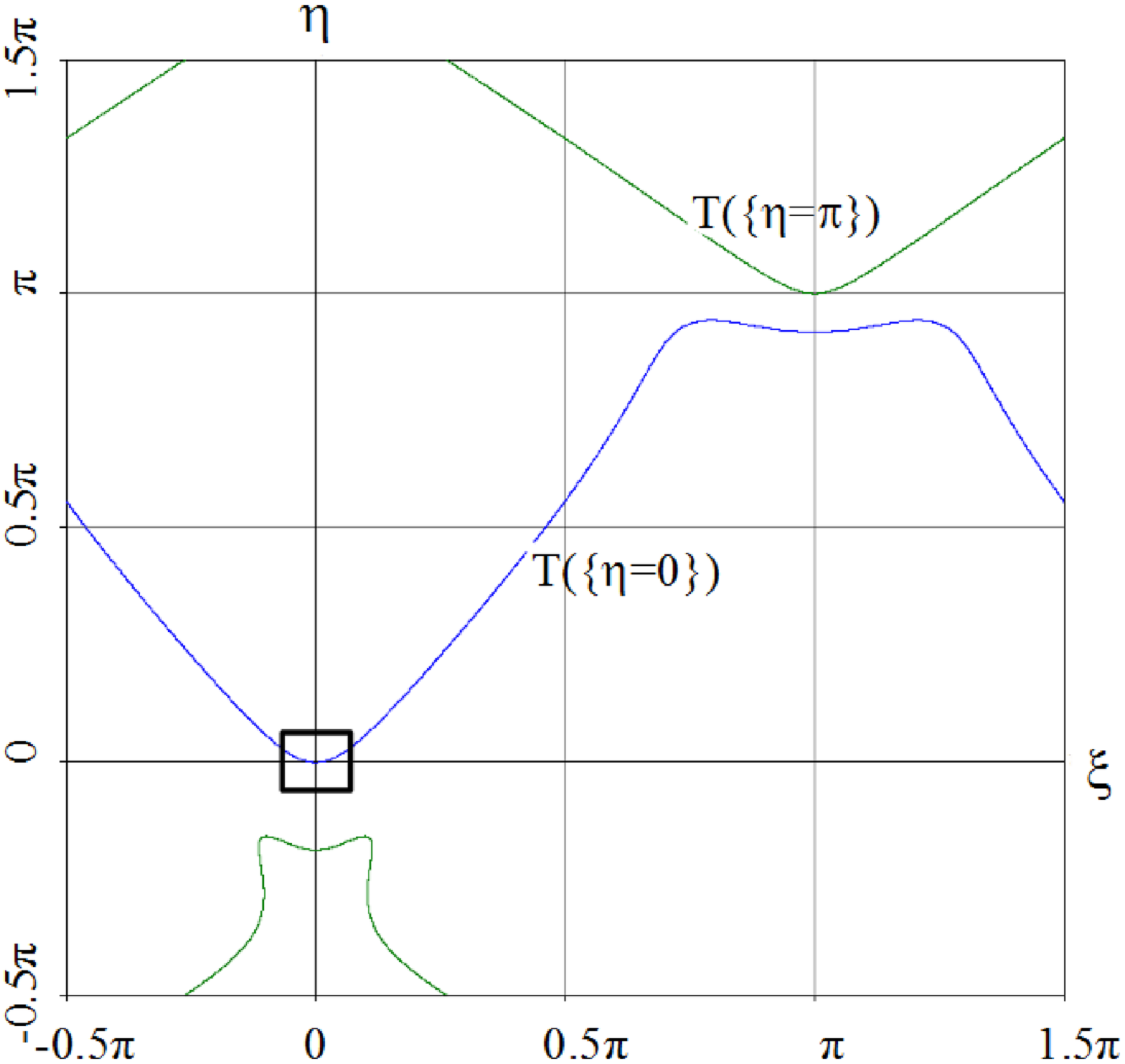} \\ (a)}
\end{minipage}
\hfill
\begin{minipage}[h]{0.48\linewidth}
\center{\includegraphics[width=1\linewidth]{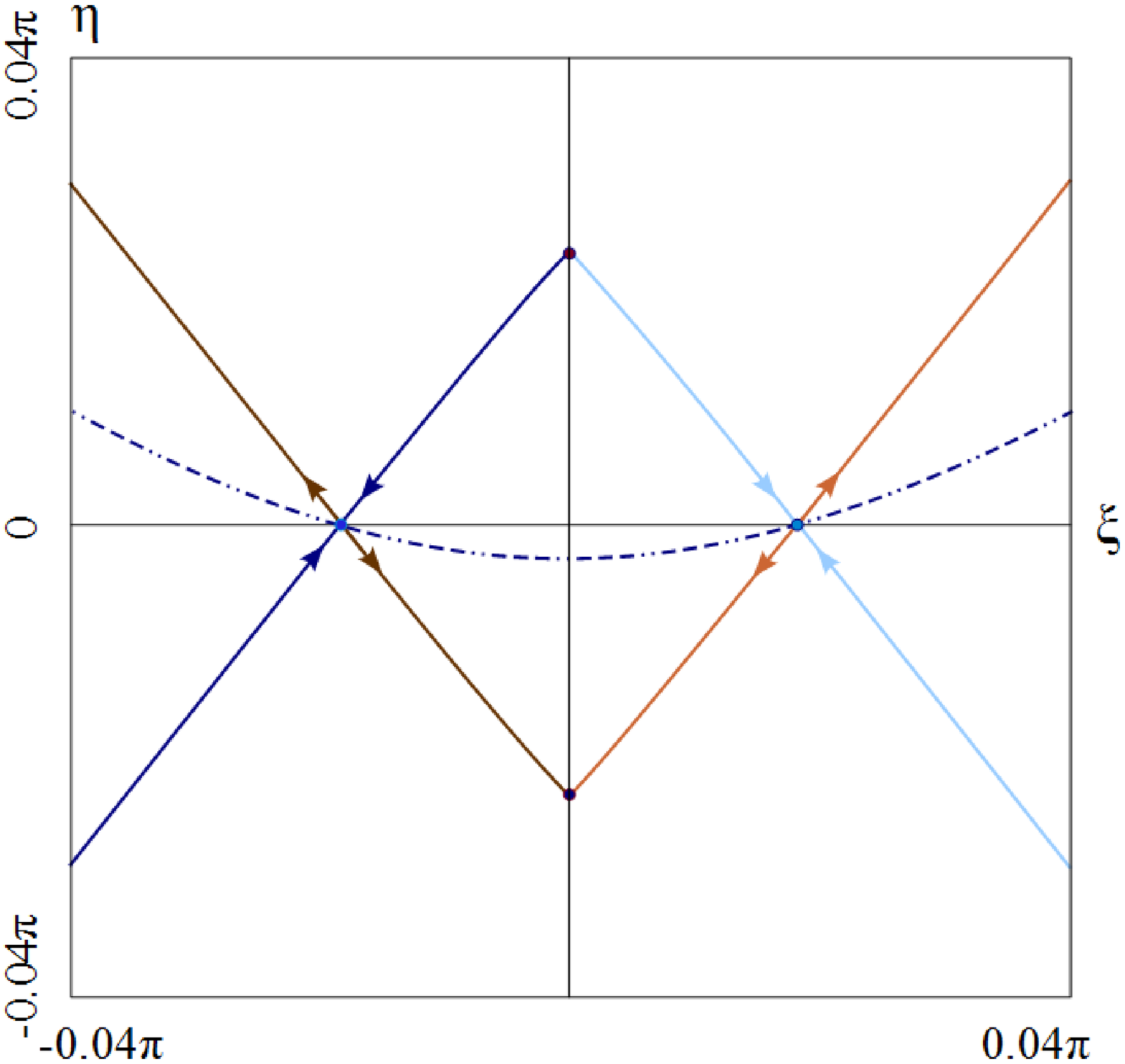} \\ (b)}
\end{minipage}
\caption{{\footnotesize a) $T(Fix(R))$ is tangent to $Fix(R)$ at two points $(\xi,\eta)=(0,0)$ and $(\xi,\eta)=(\pi,\pi)$ at
$\varepsilon = \varepsilon_1^* \simeq 0.6042$; the upper and bottom curves are the $T$-image of the lines $\eta=\pi$ and $\eta=0$, respectively.
b) The phase portrait near $(\xi,\eta)=(0,0)$ at  $\varepsilon < \varepsilon_1^*$. There are 4 fixed points, 2 symmetric saddles
on the line $\eta = 0$ and a pair of attractor and repeller on the line $\xi = 0$.}}
\label{Fig:FP1}
\end{figure}

Using this approach, we found that a pair of symmetric non-hyperbolic fixed points of the map $T$ is born at $\varepsilon = \varepsilon_1^* \simeq 0.6042$.
These points are $S$-symmetric, as they belong to the symmetry lines $\xi=\pi$ and $\xi=0$ (see Fig.~\ref{Fig:FP1}a). Recall that the map $T$
is a square of the orientation-reversing map $T_\star$. The fixed points of $T$ we discuss here form an $S$-symmetric period-2 orbit of $T_\star$.
The symmetry implies that the bifurcation of the birth of this orbit is described by the normal form (\ref{potok}). Indeed, we see that, as $\varepsilon$ increases,
both the fixed points of $T=T_\star^2$ undergo a symmetry-breaking bifurcation of type
$0\to 4$: at $\varepsilon > \varepsilon_1^*$ each point breaks into four fixed points - two $R$-symmetric (i.e. belonging to $Fix R$)
saddle points and an attractor-repeller pair on the invariant line $\xi=0$ or $\xi=\pi$, see Fig.~\ref{Fig:FP1}b.

At $\varepsilon>\varepsilon_1^*$ the system has a clearly visible sink-source pair born at $\varepsilon=\varepsilon_1^*$.
At $\varepsilon<\varepsilon_1^*$ we did not detect sinks or sources. Instead, we establish their existence in an indirect way, by finding pairs
of $R$-asymmetric saddles connected by orbits of heteroclinic tangency
(when such tangencies split as $\varepsilon$ varies sinks and sources are born \cite{DGGLS13,GST97,LSt04}). We find such pairs of saddles near
the bifurcation moments when they split off a degenerate $R$-symmetric periodic orbit.

\begin{figure}[h!]
\begin{minipage}[h]{0.48\linewidth}
\center{\includegraphics[width=1\linewidth]{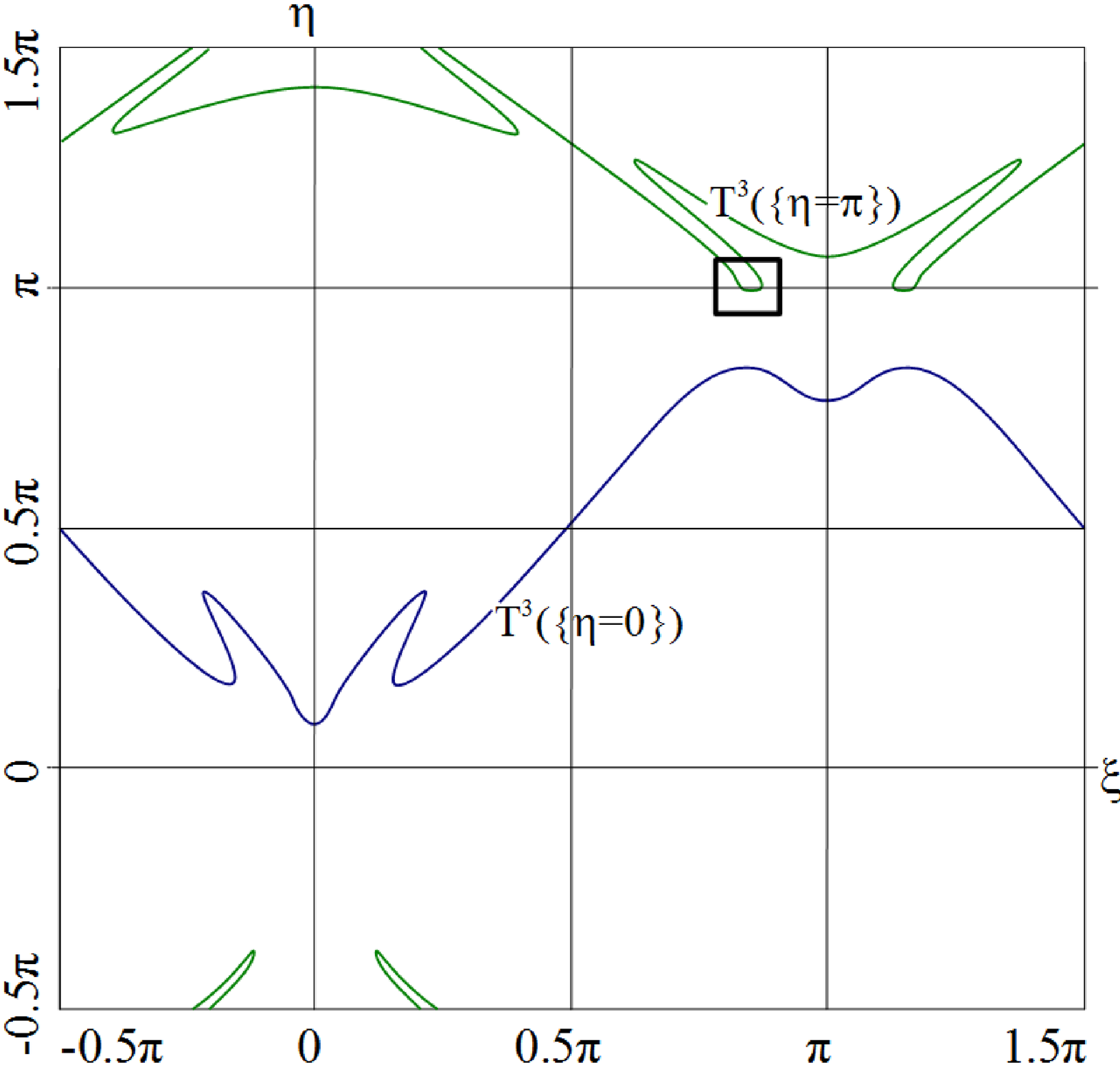} \\ (a)}
\end{minipage}
\hfill
\begin{minipage}[h]{0.48\linewidth}
\center{\includegraphics[width=1\linewidth]{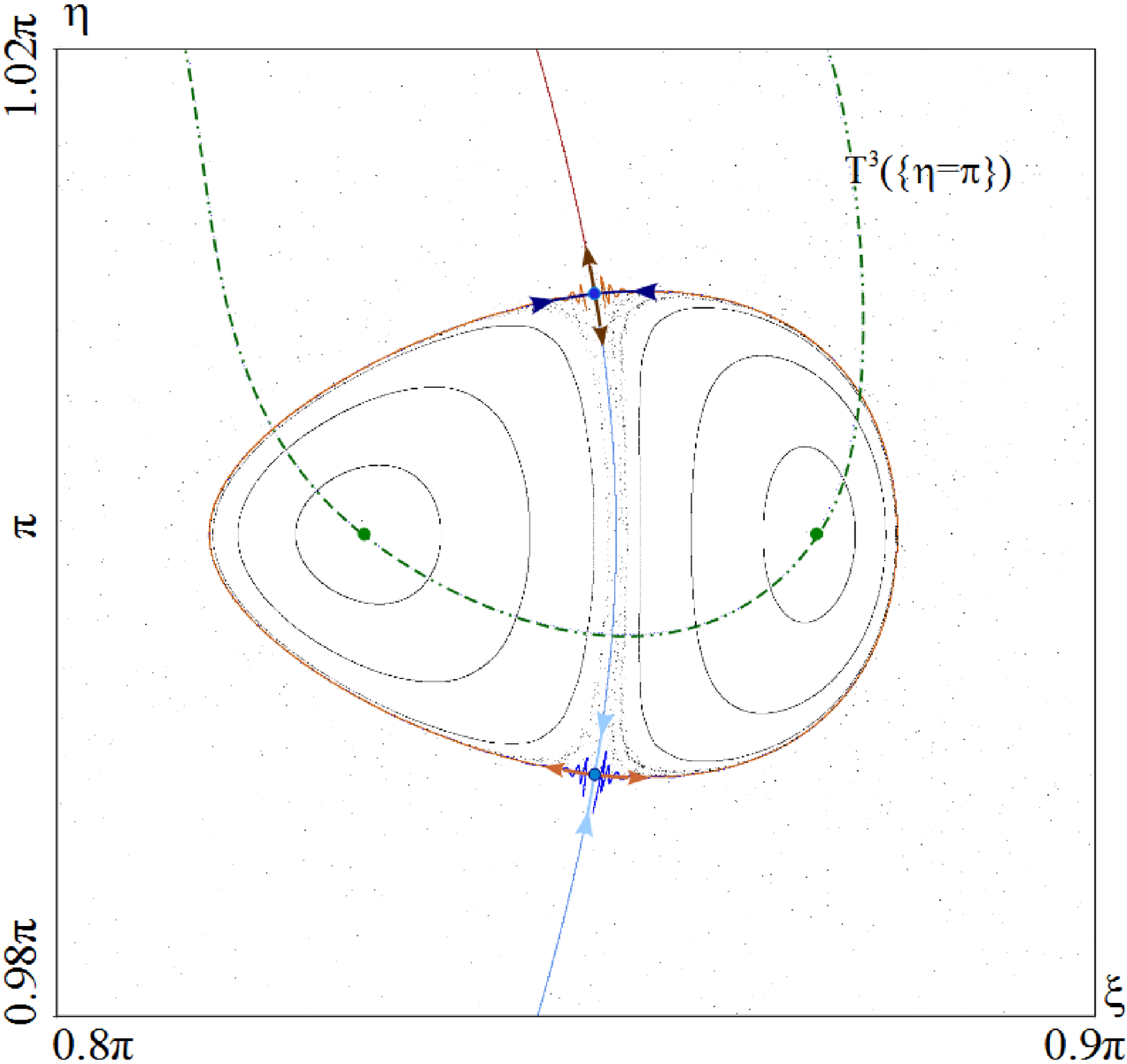} \\ (b)}
\end{minipage}
\caption{{\footnotesize
a) The image of $T^3(Fix(R))$ at $\varepsilon = 0.457 > \varepsilon_{31}^*$. Here $T^3(Fix(R))$ intersects the line $\eta = \pi$ at four points,
pairwise symmetric with respect to the line $\xi = \pi$.
b) A magnification of the rectangle region from Fig.~a. In this region there are 4 fixed points for $T^3$, two symmetrical elliptic ones and a pair of
non-conservative saddles.
}}
\label{Fig:FP3_1}
\end{figure}

At $\varepsilon = \varepsilon_{31}^* \simeq 0.455$ we found the bifurcation of the birth of an $R$-symmetric orbit of period 3 (see Fig.~\ref{Fig:FP3_1}).
At this value of $\varepsilon$ the curve $T^3(\{\eta = \pi\})$ touches the line $\eta = \pi$ at two points symmetric
with respect to the line $\xi = \pi$. We checked
that these points corresponds to different (symmetric to each other by $S$) orbits of period 3 for the map $T_{\star}$, so they are fixed points of
$T_\star^3$. Since the map $T_{\star}$ is orientation-reversing, the multipliers of these fixed points are $(+1,-1)$. As $\varepsilon$ increases, each
point splits into 4 fixed points for the map $T^3$: 2 symmetrical elliptic points and a symmetric pair of saddles (see Fig.~\ref{Fig:FP3_1}b).
For the map $T_{\star}^3$ we have, respectively, two saddle fixed points and one elliptic cycle of period two, i.e., this bifurcation
corresponds to $\alpha<0$, see Section~\ref{sec:bnorm}. Note that for $\varepsilon= 0.457$, as in Fig.~\ref{Fig:FP3_1},
the Jacobian of the derivative of the first-return map $T^3$ at the upper saddle is $\simeq 0.7198$ and the Jacobian at the lower saddle is $\simeq 1.3893$.
Thus, we indeed have here a symmetry-breaking bifurcation, which creates an $R$-symmetric pair of saddles with Jacobians different from $1$, i.e. breaks
the conservativity and creates conditions for the birth of long-period sinks and sources via bifurcations of heteroclinic and homoclinic tangencies to the saddles.

\begin{figure}[h!]
\begin{minipage}[h]{0.48\linewidth}
\center{\includegraphics[width=1\linewidth]{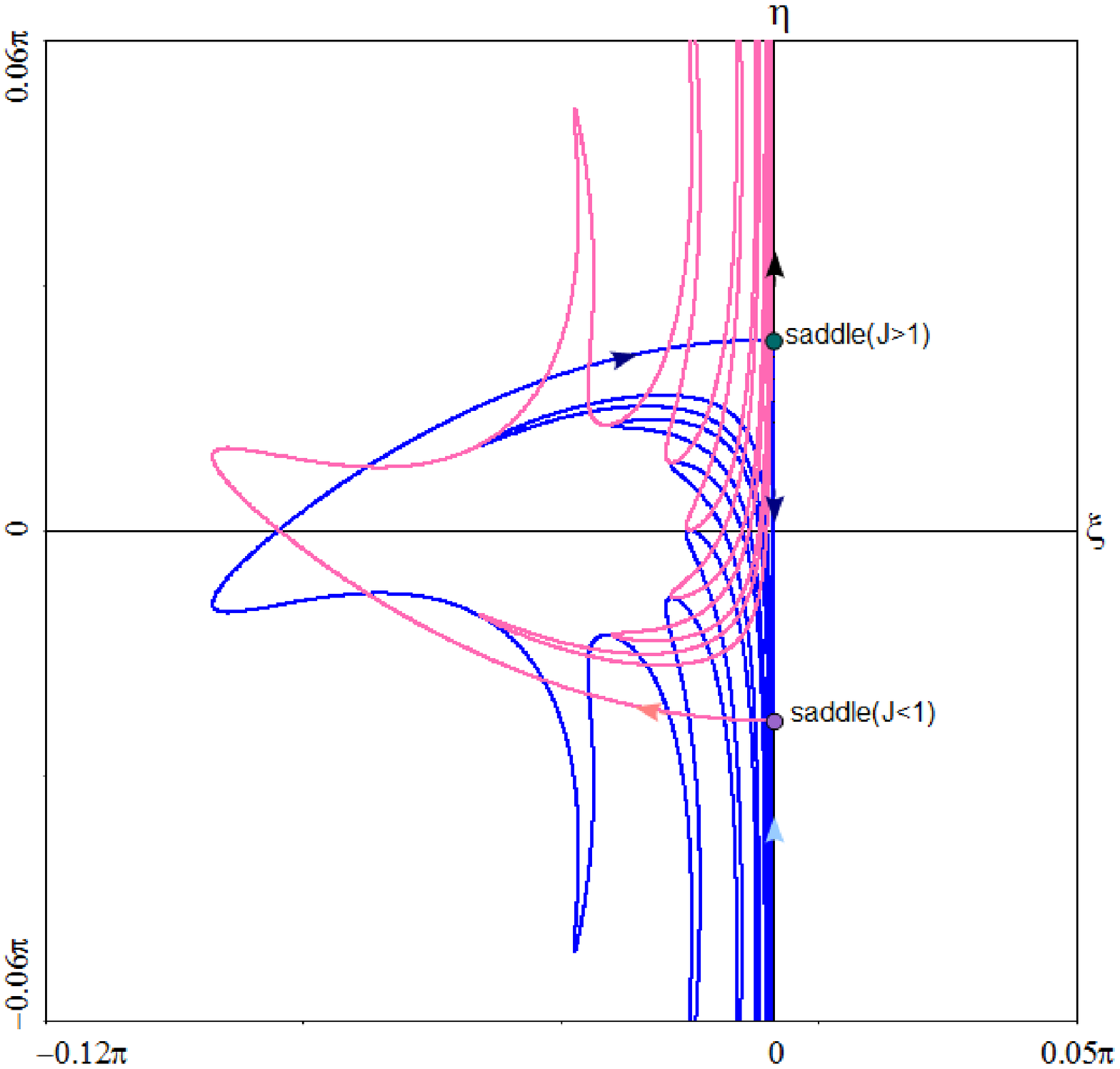} \\ (a)}
\end{minipage}
\hfill
\begin{minipage}[h]{0.48\linewidth}
\center{\includegraphics[width=1\linewidth]{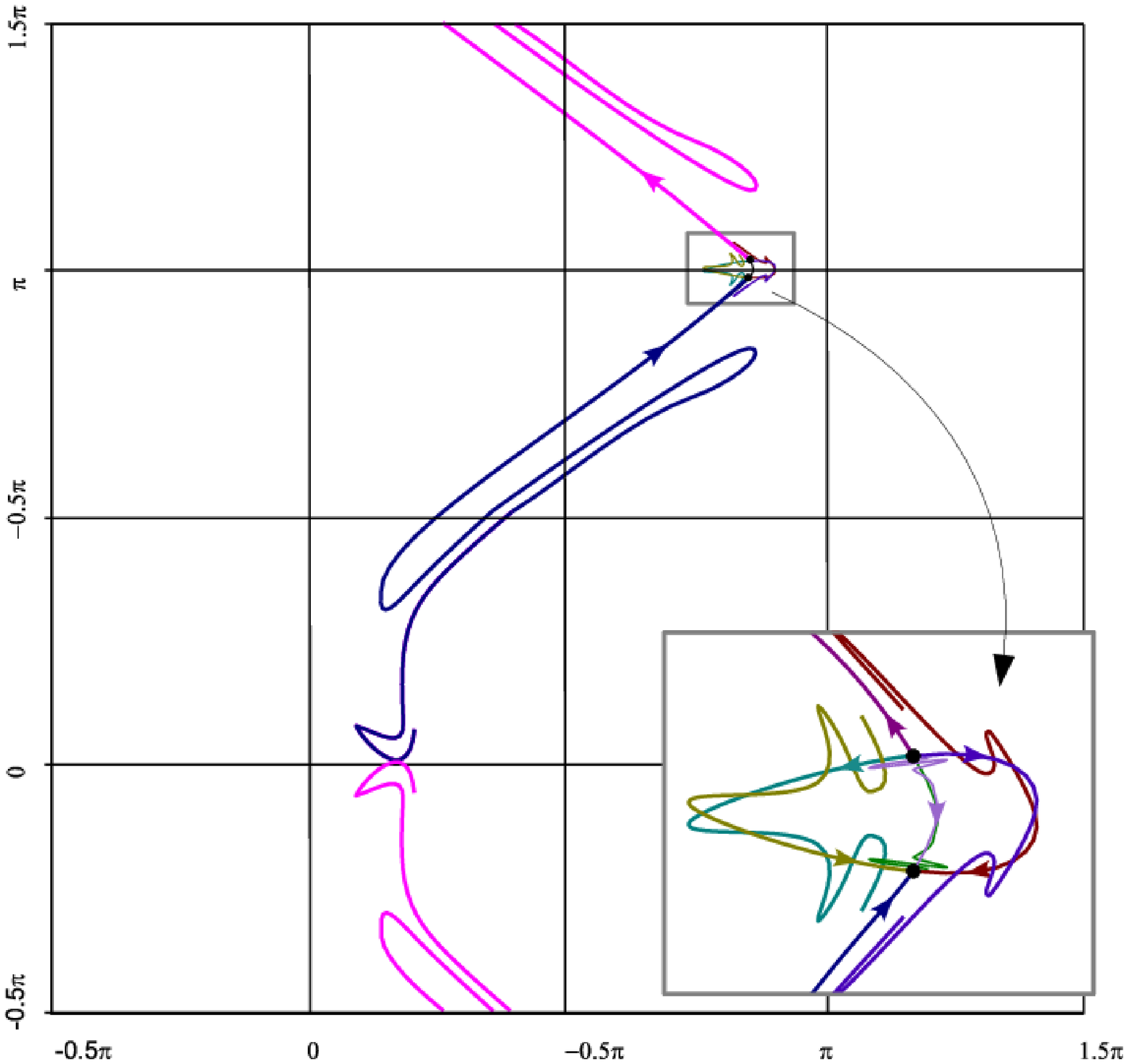} \\ (b)}
\end{minipage}
\caption{{\footnotesize
Creating of ``large'' heteroclinic orbits at $\varepsilon \approx 0.463$.}}
\label{Fig:Heteroclinic_large}
\end{figure}

There are two ways such tangencies can appear near this bifurcation. The first way gives ``small'' heteroclinic orbits which emerge due to the
splitting of the separatrix connection which exists in the flow normal form (\ref{potok}) as in Fig.~\ref{4typlocbif}a. We found examples of these orbits
at $\varepsilon\simeq 0.487$, see Fig.~\ref{Fig:Heteroclinic_large}a. Another way leads to creation of ``large'' heteroclinic orbits. At the moment of bifurcation the degenerate
periodic point has one unstable separatrix and one stable separatrix which leave a small neighborhood of the point. If there is
no additional symmetry, there is no reason why these separatrices cannot have intersections. After the bifurcation, these separatrices become an unstable separatrix
of one saddle and a stable separatrix of the other saddle; they may have transverse intersections and tangencies outside a small neighborhood of the saddles.
We find an example of such tangency at $\varepsilon\simeq 0.463$, see Fig.\ref{Fig:Heteroclinic_large}b. The corresponding values of Jacobians are $J_1\simeq 0.524<1$, $J_2\simeq 1.909>1$.
In both these examples, the heteroclinic tangencies are constituent parts of heteroclinic cycles that
include saddle points with Jacobians less and larger than $1$, so splitting the tangencies as $\varepsilon$
varies must lead to the birth of stable and unstable periodic orbits.

\begin{figure}[h]
\begin{minipage}[h]{0.47\linewidth}
\center{\includegraphics[width=1\linewidth]{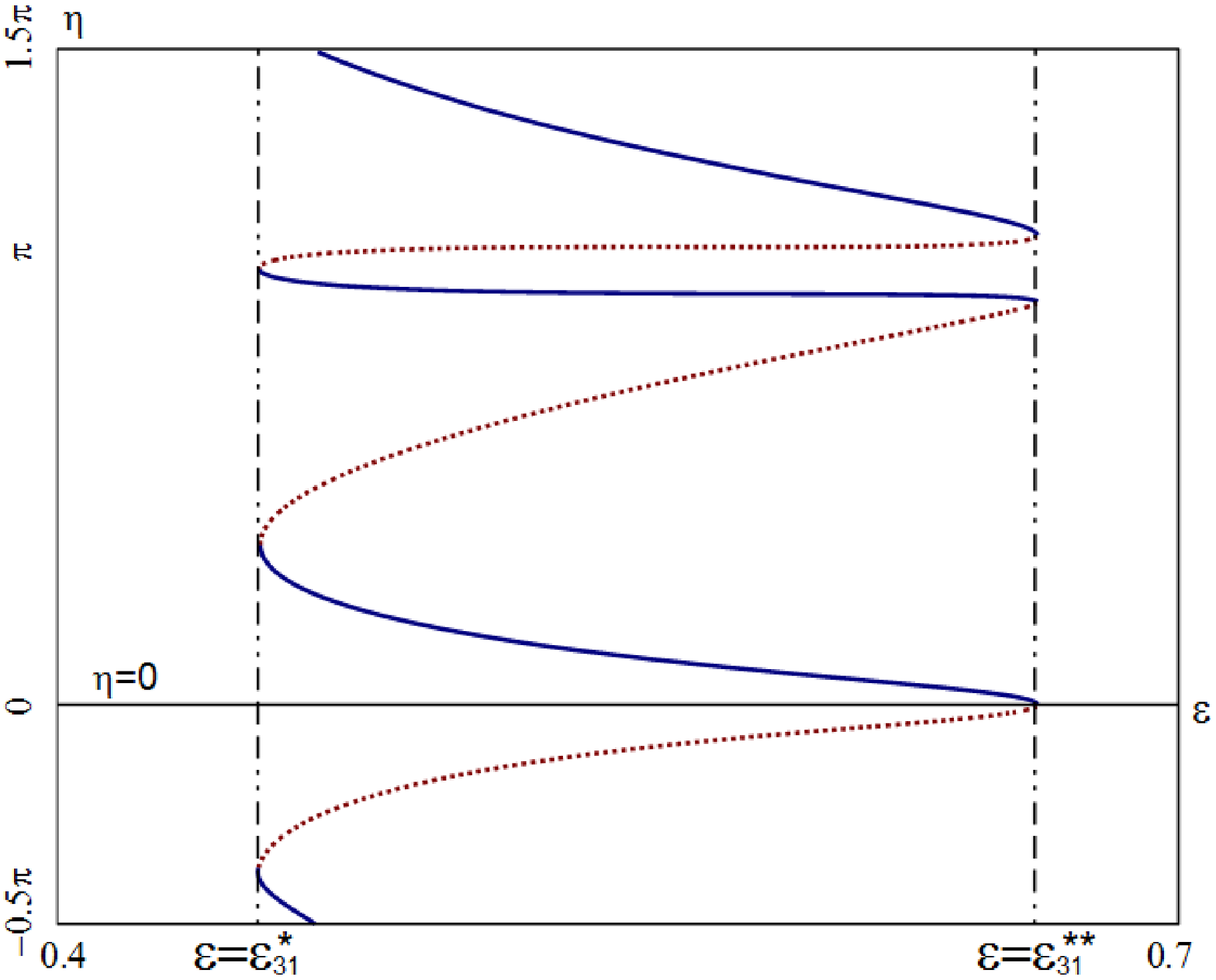} (a)}
\end{minipage}
\hfill
\begin{minipage}[h]{0.49\linewidth}
\center{\includegraphics[width=1\linewidth]{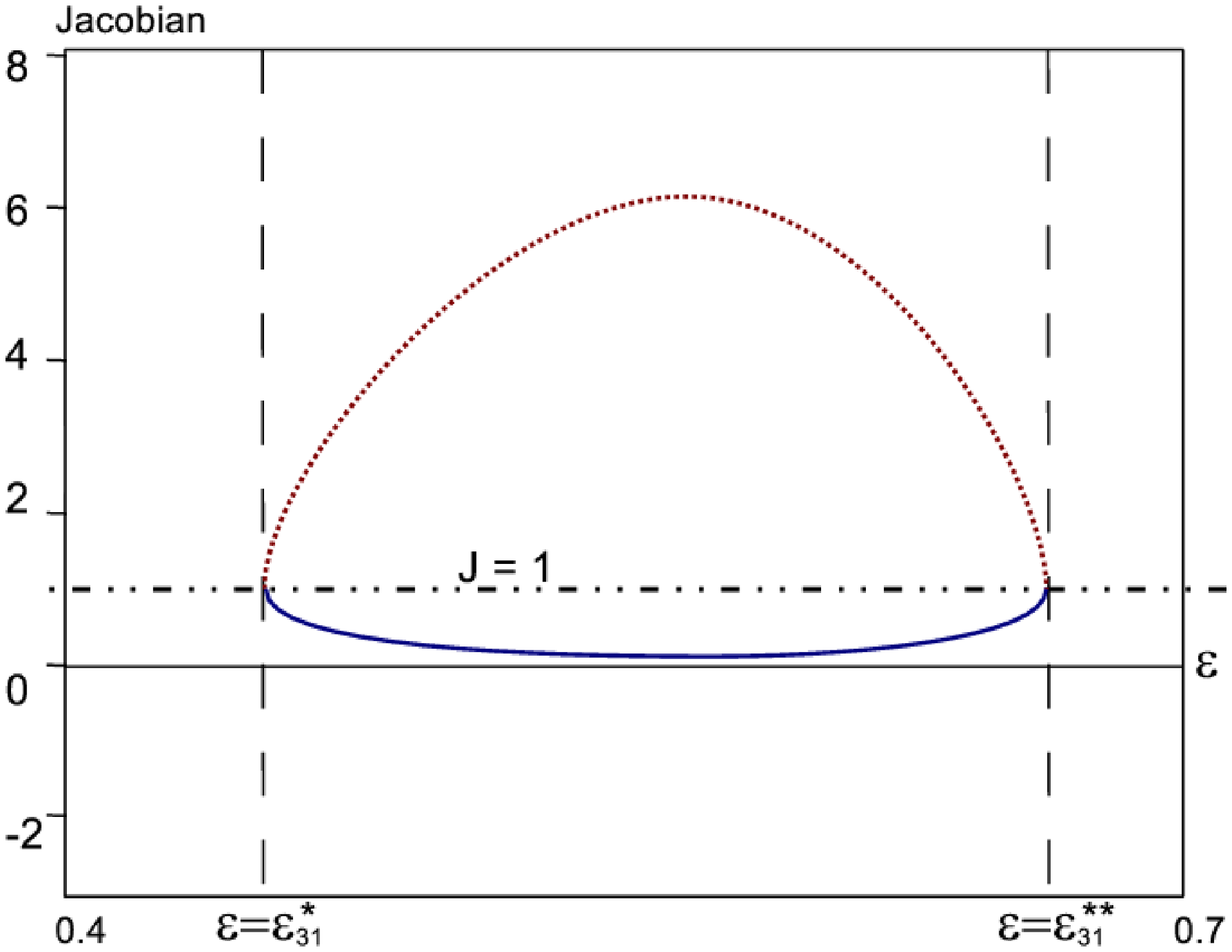} (b)}
\end{minipage}
\caption{{\footnotesize (a) The evolution of the pair of non-conservative period-3 saddles from varying of $\varepsilon$ from  $\varepsilon_{31}^* \approx 0.457$ to $\varepsilon_{32}^* \approx 0.664$; (b) graphs of the Jacobians for these saddles.}}
\label{Fig:BiffDiagrAndJacobian}
\end{figure}

We traced the evolution of the pair of non-conservative period-3 saddles as the parameter $\varepsilon$ increases. We found, see Fig.~\ref{Fig:BiffDiagrAndJacobian}, that the saddles exist until
$\varepsilon=\varepsilon_{30}^*\simeq 0.663$ when they collide with each other (and a pair of $R$-symmetric elliptic period-3 points) at a $0\to 4$
bifurcation and disappear. For all $\varepsilon\in(\varepsilon_{31}^*,\varepsilon_{30}^*)$ the Jacobians of the saddles remain different from $1$, and the
stable and unstable manifolds of the saddles have intersections. Tangencies between the invariant manifolds of the saddles appear easily for systems without
an uniformly hyperbolic structure, so we may conjecture that non-transverse heteroclinic cycles involving these saddles exist for a dense subset of the interval
$\varepsilon\in(\varepsilon_{31}^*,\varepsilon_{30}^*)$, hence the mixed dynamics exists for all $\varepsilon$ from this interval.

\begin{figure}[h!]
\begin{minipage}[h]{0.48\linewidth}
\center{\includegraphics[width=1\linewidth]{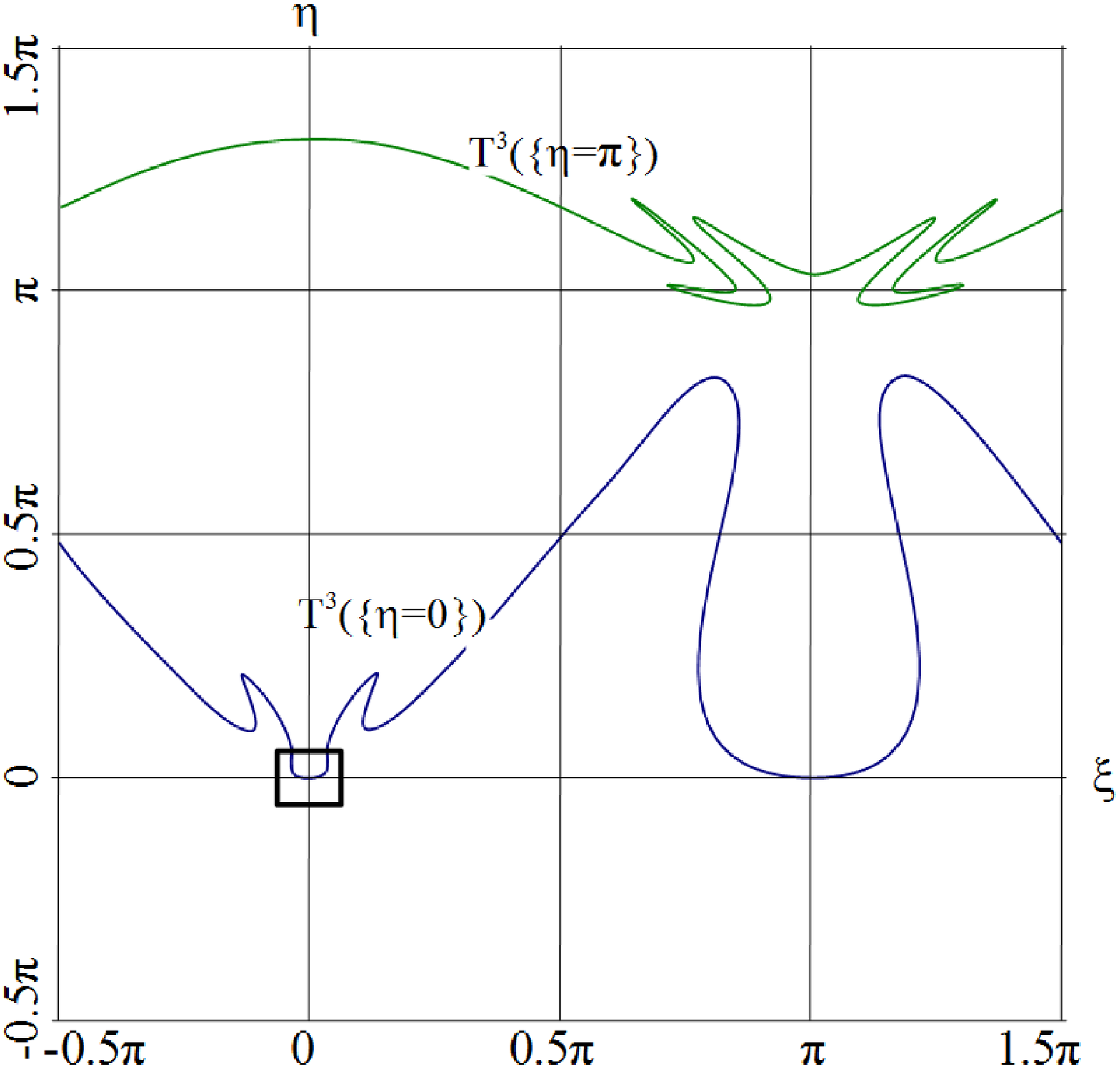} \\ (a)}
\end{minipage}
\hfill
\begin{minipage}[h]{0.48\linewidth}
\center{\includegraphics[width=1\linewidth]{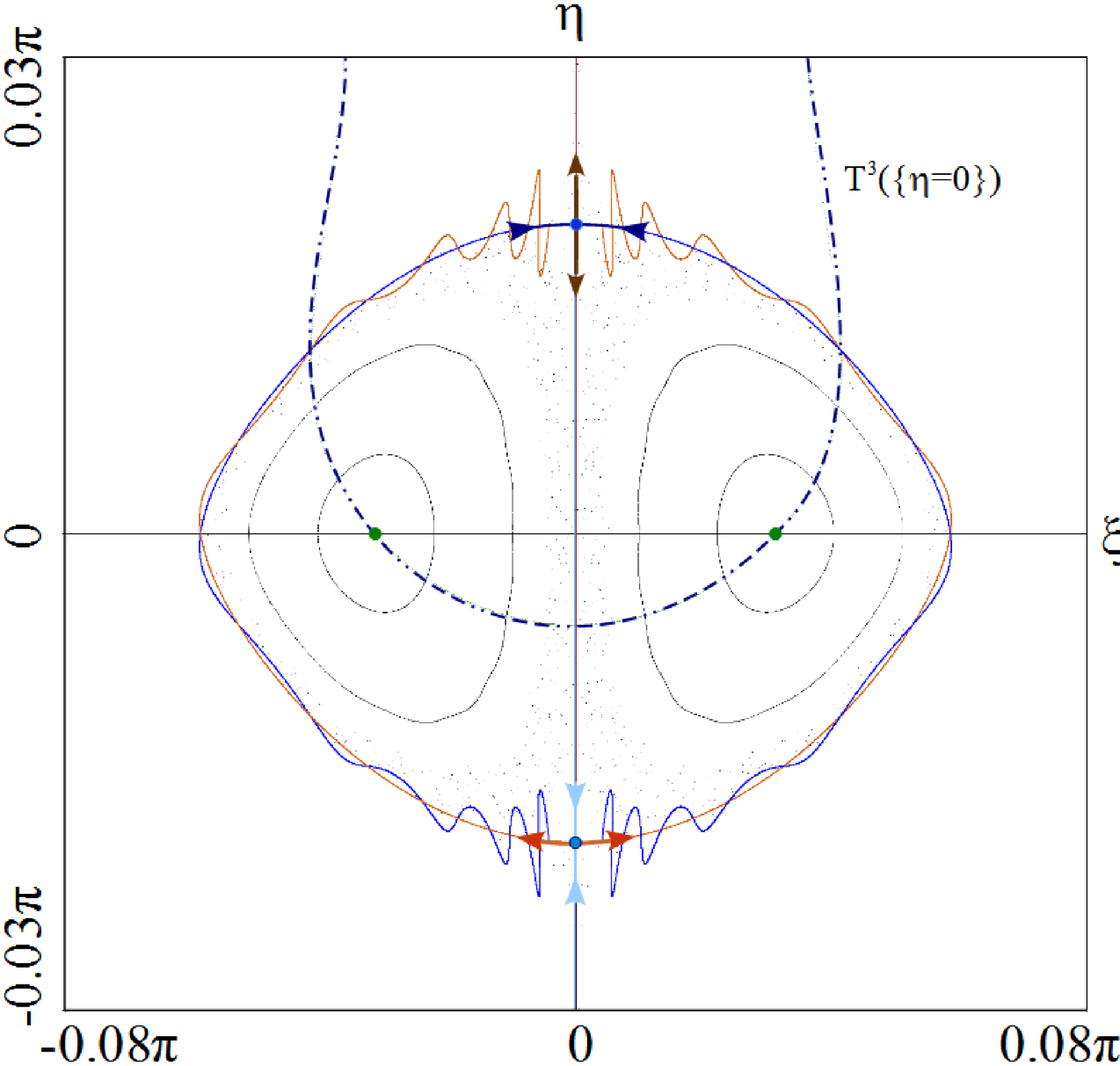} \\ (b)}
\end{minipage}
\caption{{\footnotesize
The same as the previous figure for $\varepsilon = 0.483 > \varepsilon_{32}^*$.}}
\label{Fig:FP3_2}
\end{figure}

\begin{figure}[h!]
\begin{minipage}[h]{0.48\linewidth}
\center{\includegraphics[width=1\linewidth]{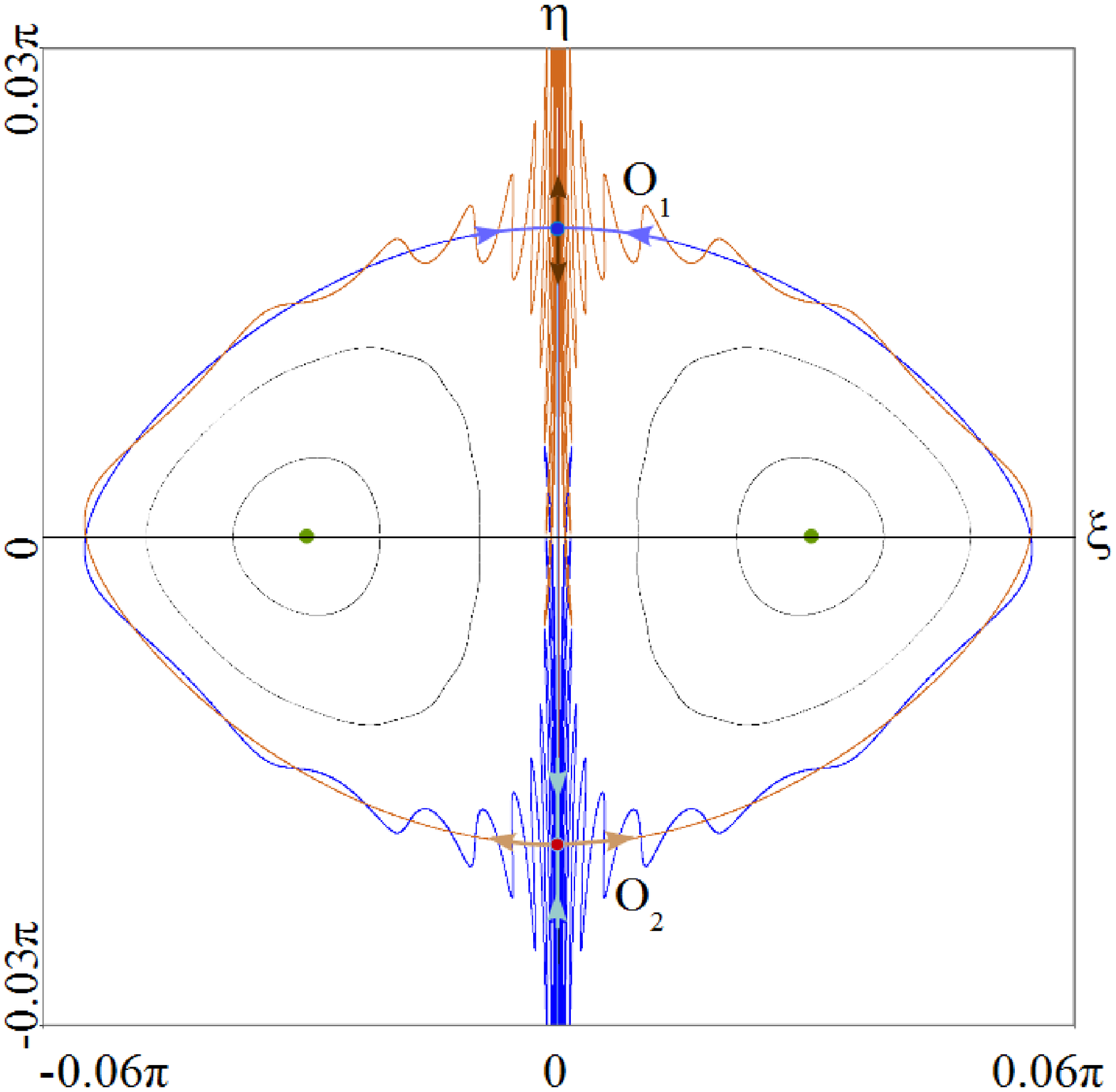} \\ (a)}
\end{minipage}
\hfill
\begin{minipage}[h]{0.48\linewidth}
\center{\includegraphics[width=1\linewidth]{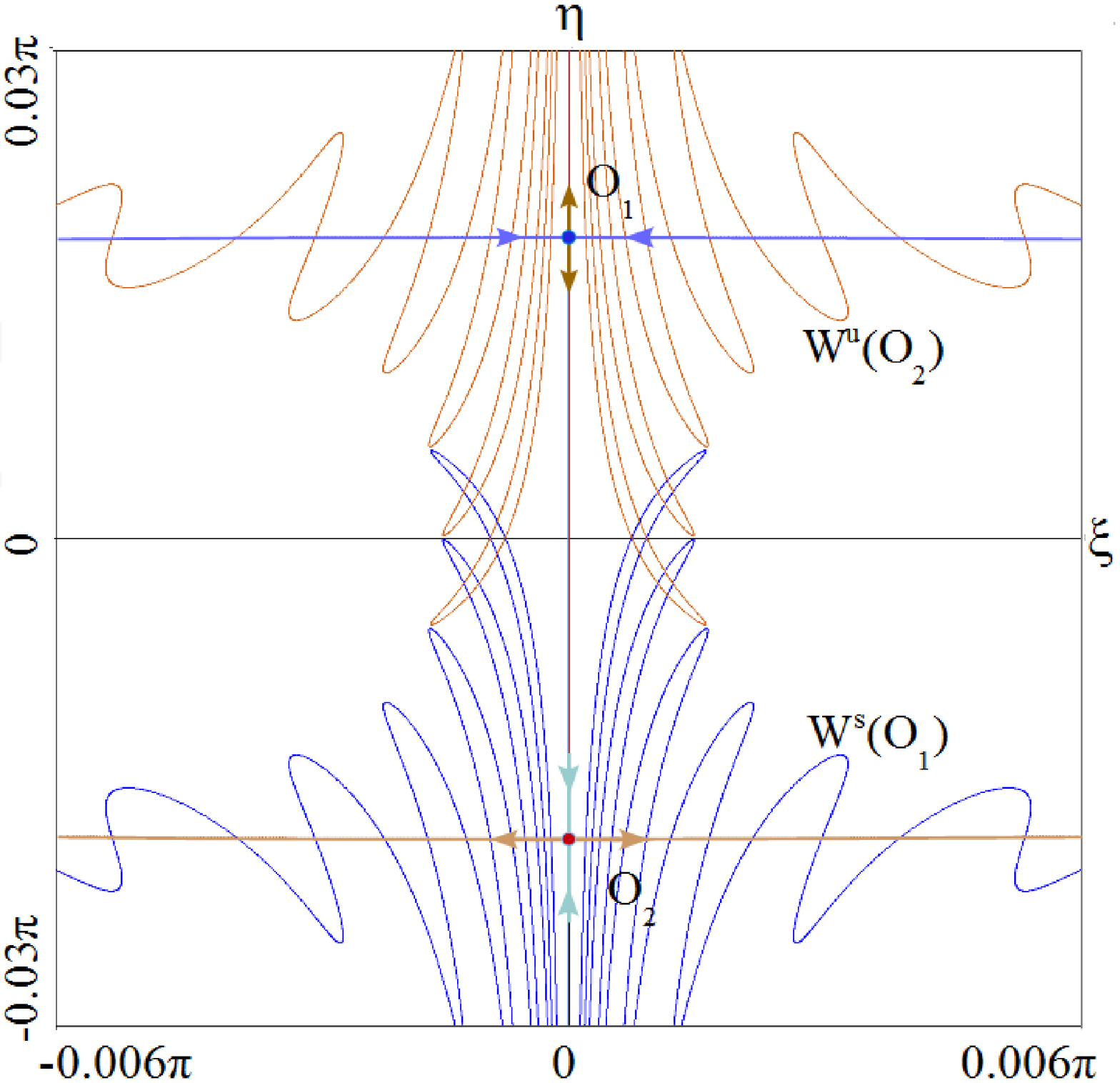} \\ (b)}
\end{minipage}
\caption{{\footnotesize a) The stable and unstable
manifolds of saddle fixed points of $T^3$ are shown at $\varepsilon = 0.485$ which compose a ``heteroclinic tangle''.  b) A magnification of fig.~(a), it is seen that a nontransversal heteroclinic cycle (of Lamb-Stenkin type) is created.}}
\label{Fig:FP3_heteroclinic}
\end{figure}

One more bifurcation of the birth of period-3 point is detected at $\varepsilon = \varepsilon_{32}^* \simeq 0.483$. Unlike the previous case,
here the curve $T^3(\{\eta = 0\})$ touches the line $\eta = 0$ at two $S$-symmetric points,
($\xi = 0, \eta = 0$) and ($\xi = \pi, \eta = 0$), see Fig.~\ref{Fig:FP3_2}a. These points form an $S$-symmetric period-2 orbit of the map $T_{\star}^3$.
Due to the symmetry $S$, this orbit splits into four after the bifurcation. So, for the map $T$,
we see that each of the period-3 orbits splits into 4 period-3 orbits, two $R$-symmetric elliptic orbits, and an $R$-symmetric pair of saddles,
see Fig.~\ref{Fig:FP3_2}b. We checked that the saddles are non-conservative, with the Jacobians $J_1 \simeq 0.9988$ (for the upper saddle)
and $J_2 \simeq 1.0012$) (for the lower saddle) at $\varepsilon=0.485$. We have found orbits of heteroclinic tangencies between the separatrices of the saddles
at the same $\varepsilon$, see Fig. \ref{Fig:FP3_heteroclinic}b. However, the other separatrices of the saddles do not split (they coincide with the symmetry line), so the corresponding heteroclinic
cycles we have here are different from those considered in \cite{GST97, LSt04}. We, therefore, do not have a theorem which, like in the previous case, would guarantee that,
as the heteroclinic tangencies are split, periodic attractors and repellers are born. However, they most probably do; this question requires a further study.

\begin{figure}[h]
\begin{minipage}[h]{0.48\linewidth}
\center{\includegraphics[width=1\linewidth]{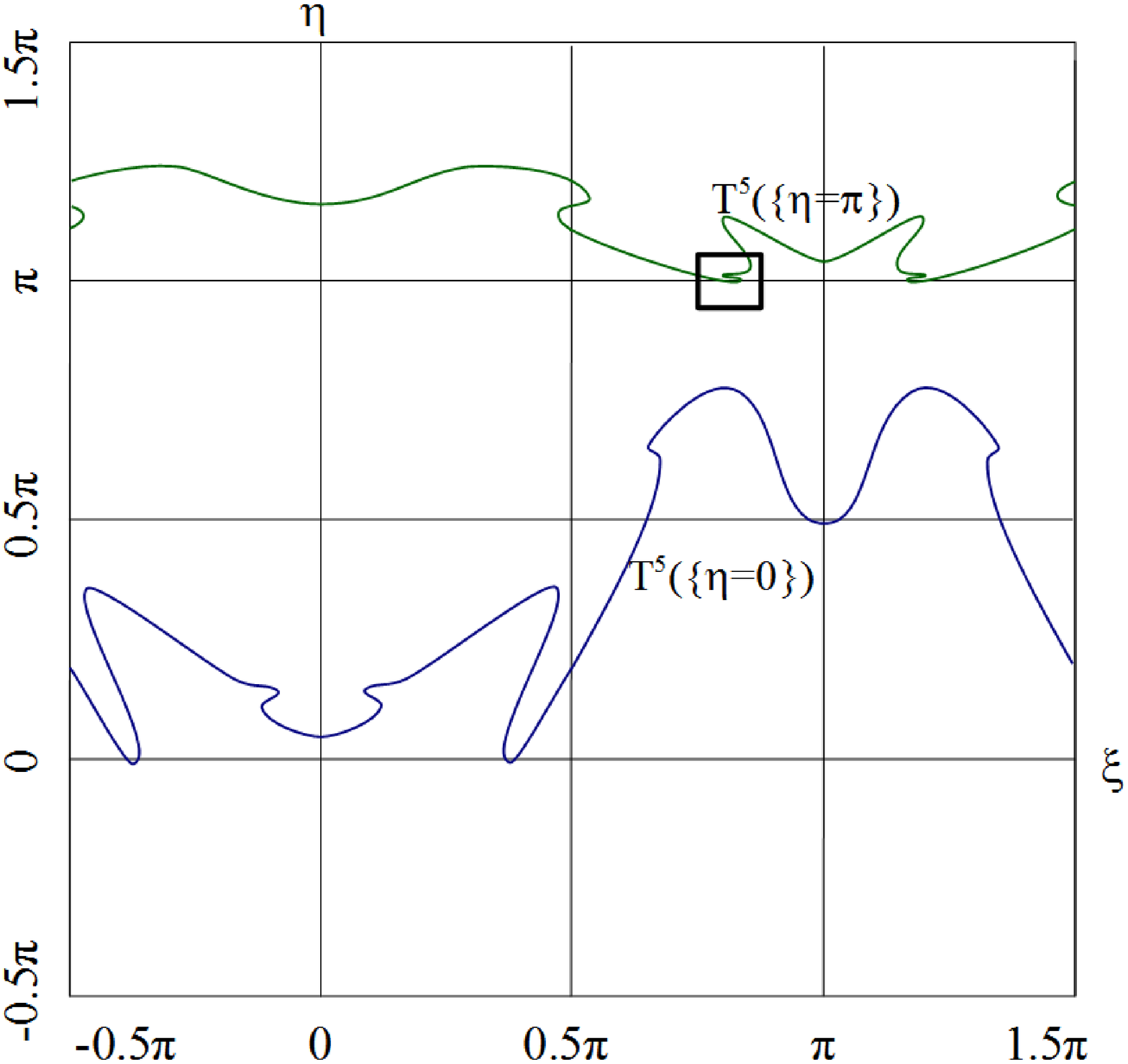} \\ (a)}
\end{minipage}
\hfill
\begin{minipage}[h]{0.48\linewidth}
\center{\includegraphics[width=1\linewidth]{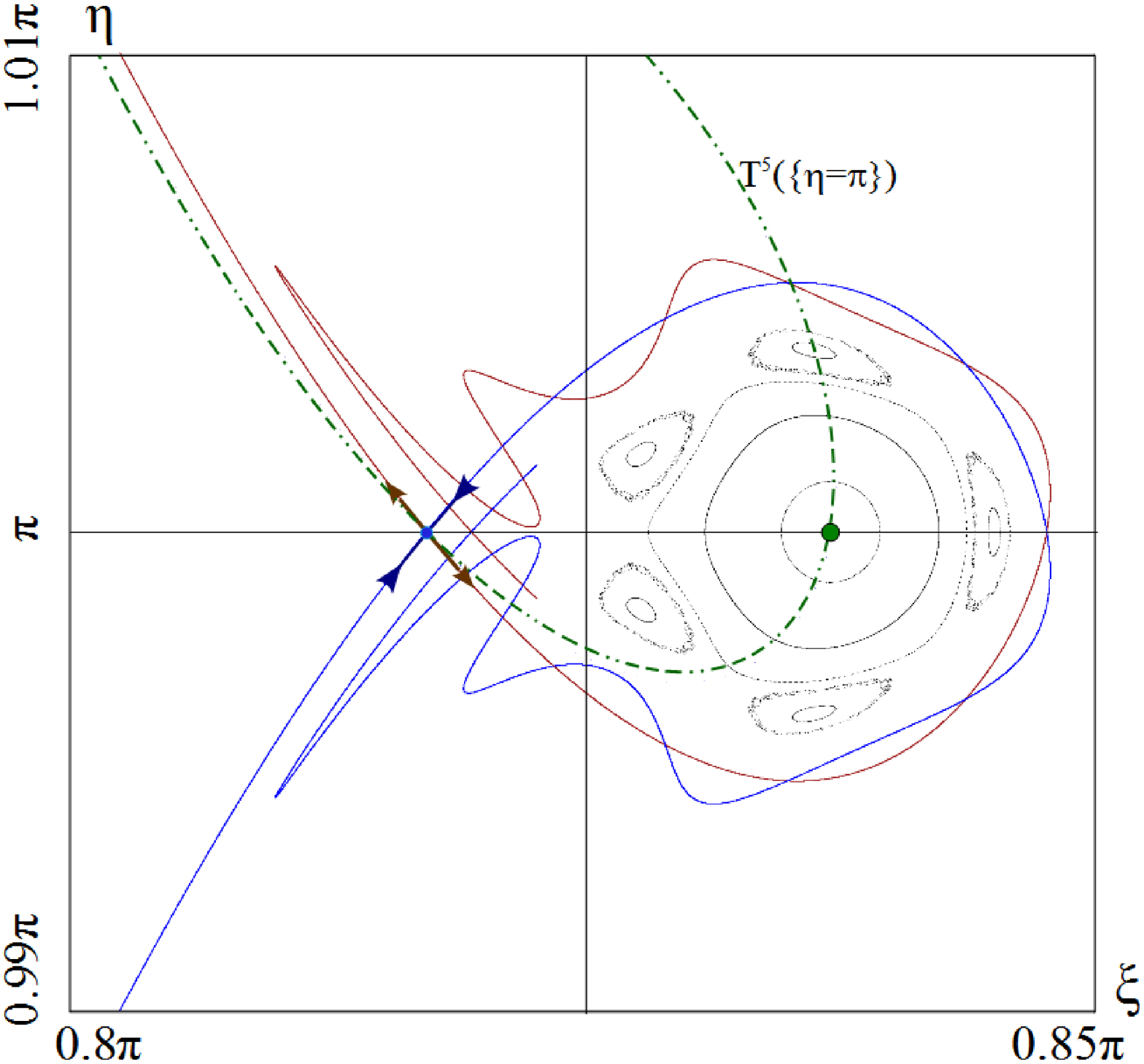} \\ (b)}
\end{minipage}
\caption{{\footnotesize
a) The image of $T^5(Fix(R))$ is shown for $\varepsilon = 0.418 > \varepsilon_{5}^*$. Here $T^5(Fix(R))$ intersects the lines $\eta = 0$ and  $\eta = \pi$ at four points each.
b) Magnification of the rectangle region from Fig.~a. In this region there are 2 fixed points for $T^5$, elliptic and saddle one, both symmetric.}}
\label{Fig:FP5}
\end{figure}

The case of period-5 points is quite different. At $\varepsilon = \varepsilon_5^* \simeq 0.417$ four tangent points
between $T^5(Fix(R))$ and $Fix(R) = \{\eta = \pi \cup \eta = 0\}$ appear (see Fig.~\ref{Fig:FP5}a). These points correspond to a pair of period-2 orbits for the
map $T_{\star}^5$; the orbits are symmetric to each other with respect to $S$. Each of these orbits is not $S$-symmetric by itself. Therefore, as $\varepsilon$ increases, we
have the usual bifurcation of a conservative-like parabolic orbit: each of the orbits splits into two symmetric orbits, one elliptic and one saddle (see Fig.~\ref{Fig:FP5}b).

\begin{figure}[ht!]
\begin{minipage}[h]{0.48\linewidth}
\center{\includegraphics[width=1\linewidth]{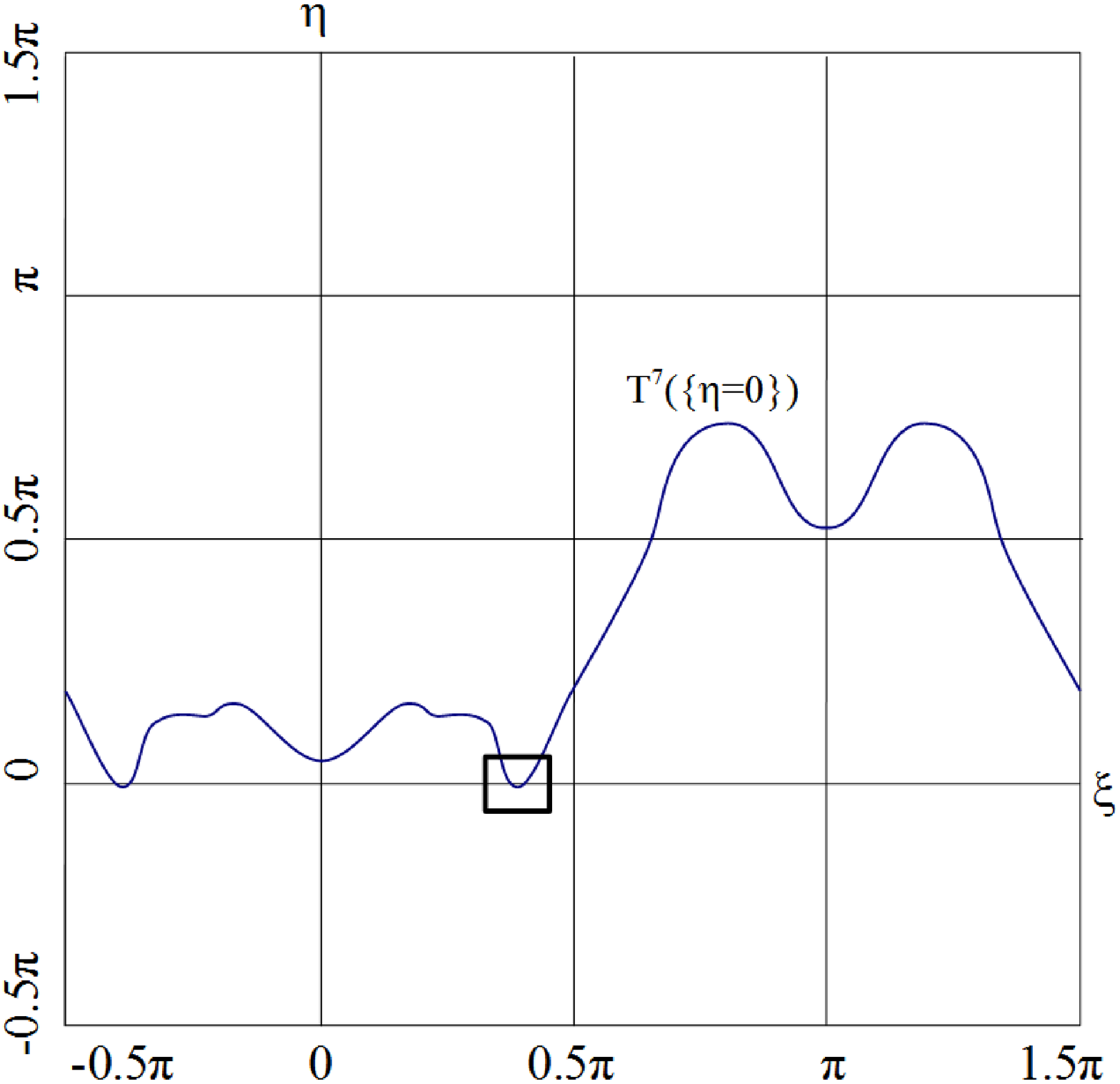} \\ (a)}
\end{minipage}
\hfill
\begin{minipage}[h]{0.48\linewidth}
\center{\includegraphics[width=1\linewidth]{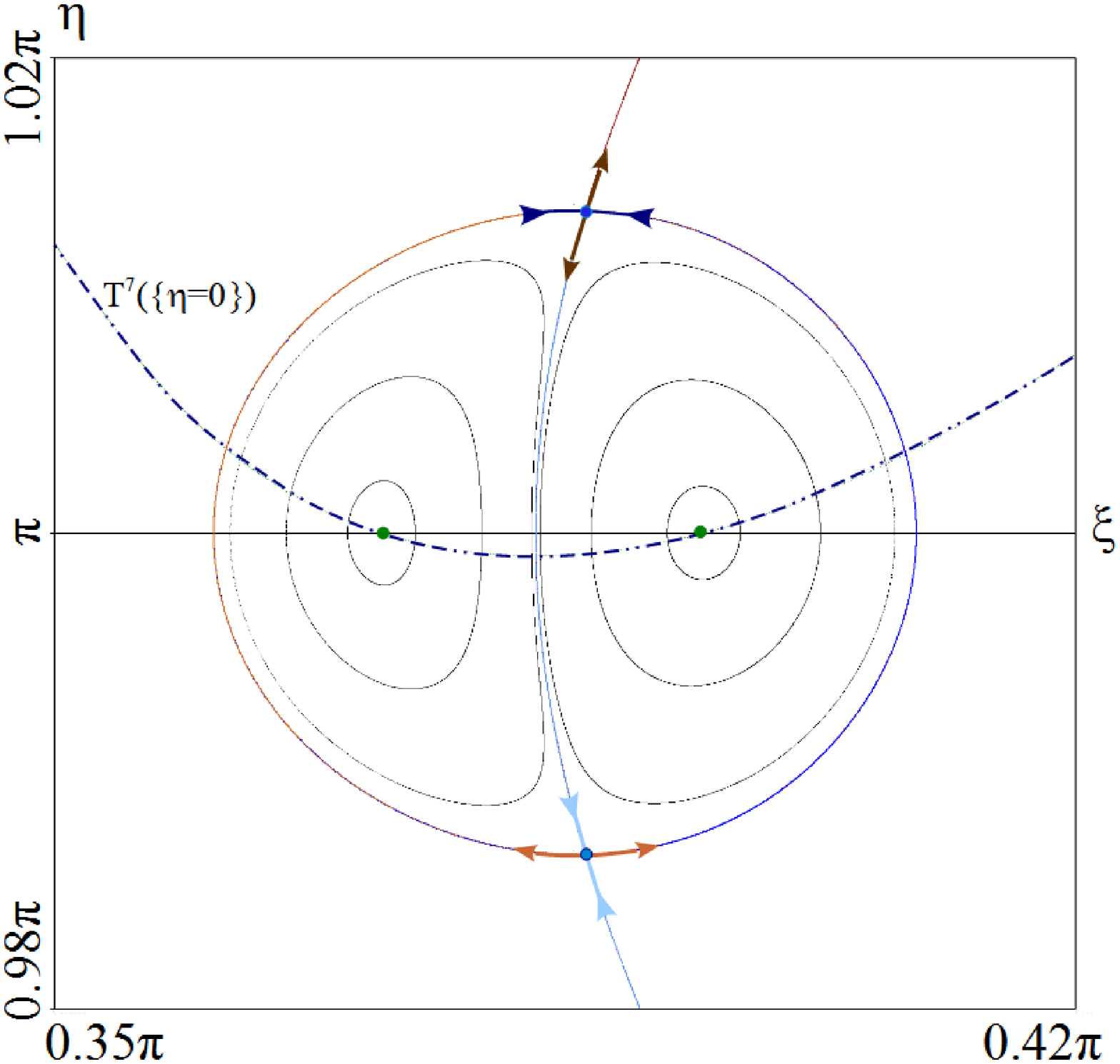} \\ (b)}
\end{minipage}
\vfill
\begin{minipage}[h]{0.48\linewidth}
\center{\includegraphics[width=1\linewidth]{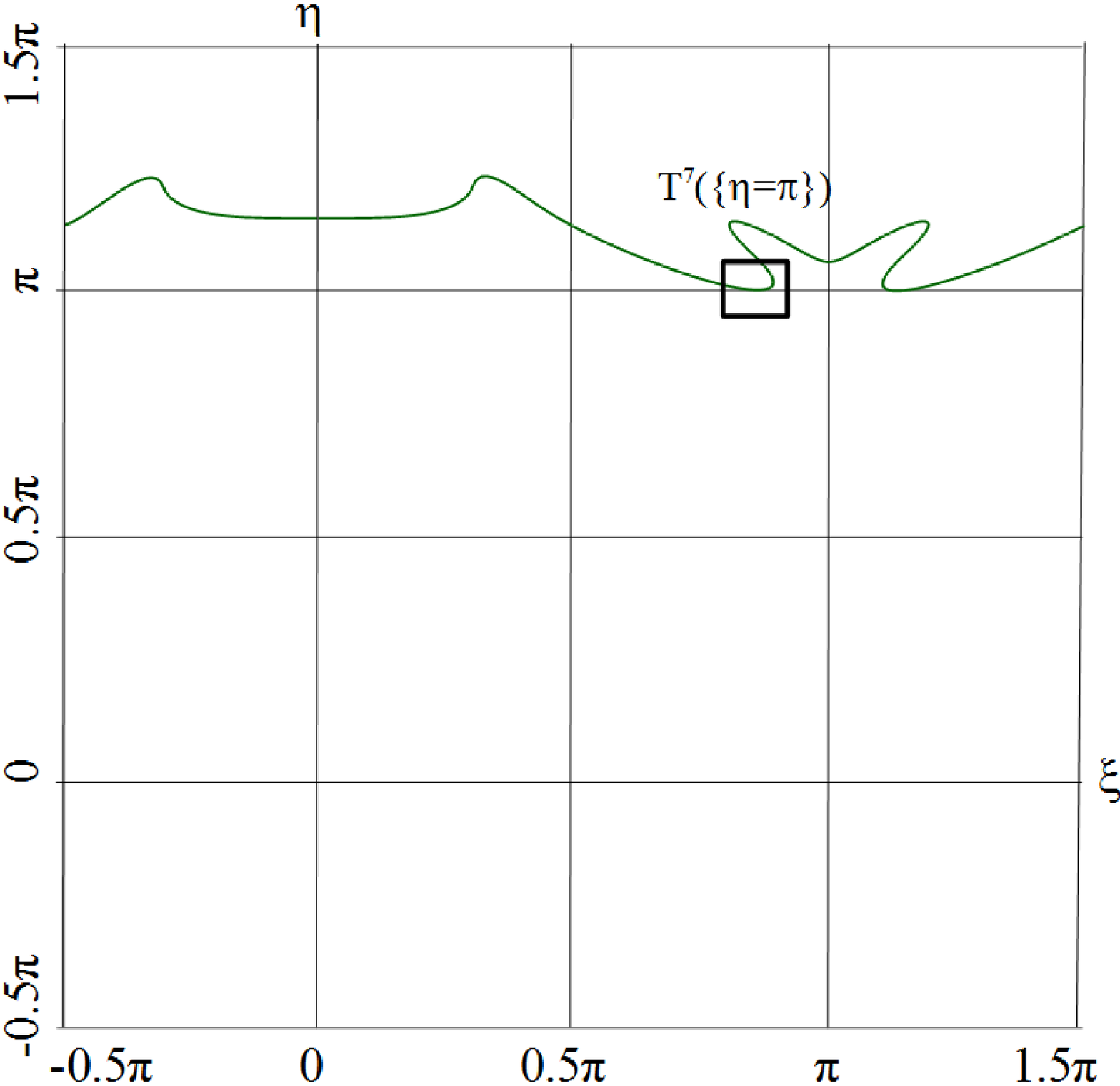} \\ (a)}
\end{minipage}
\hfill
\begin{minipage}[h]{0.48\linewidth}
\center{\includegraphics[width=1\linewidth]{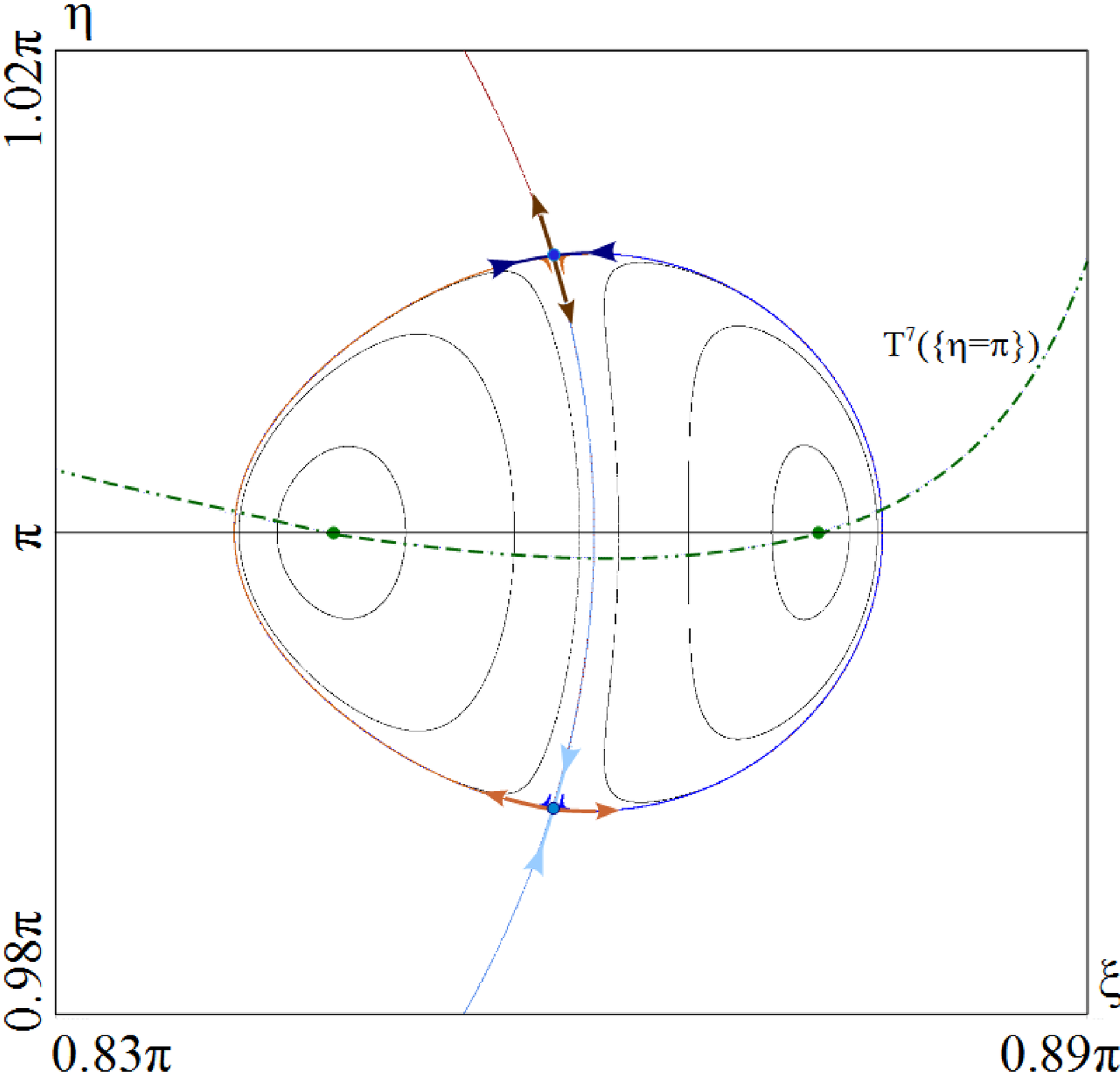} \\ (b)}
\end{minipage}
\caption{{\footnotesize
a) The image of $T^7(\{\eta = 0\})$ at $\varepsilon = \varepsilon_{71}^* \simeq 0.380$.
b) A magnification of the rectangle region from Fig.~a.
c) The image of $T^7(\{\eta = \pi\})$ at $\varepsilon = \varepsilon_{72}^* \simeq 0.381$.
b) A magnification of the rectangle region from Fig.~c.
}}
\label{Fig:FP7}
\end{figure}

Symmetric orbits of period 7 are born at $\varepsilon = \varepsilon_{71}^* \simeq 0.3795$ and $\varepsilon = \varepsilon_{72}^* \simeq 0.3805$,
at the moments of tangency of $T^7(\eta = \pi)$ with $\eta = \pi$ and $T^7(\eta = 0)$ with $\eta = 0$, respectively. Every point of these orbits
is a fixed point for $T_{\star}^7$. Therefore, with the increase of $\varepsilon$, each orbit splits into 4 fixed points of the map $T^7$,
two $R$-symmetric elliptic points, and a symmetric pair of non-symmetric saddles (this corresponds to two saddle fixed points and one elliptic period-2
orbit for $T_{\star}^7$), see Fig.~\ref{Fig:FP7}. Thus we have here the same situation as in the first case of period-3 orbits.

We have also found that bifurcations of the birth of $R$-symmetric orbits of period 9 (at $\varepsilon = \varepsilon_{9}^* \simeq 0.348$) have
the type $0\to 2$ (the same as for period 5). The birth of an orbit of period 11 was detected at $\varepsilon = \varepsilon_{11}^* \simeq 0.323$,
however we could not determine the type of the corresponding bifurcation, since the elliptic zones near this orbit happen to be extremely narrow,
of the order of numerical accuracy. We conjecture that conservativity breaking bifurcations exist in this model for arbitrarily small $\varepsilon$, i.e., the mixed dynamics
exist within apparently conservative ``chaotic sea'' for all $\varepsilon>0$.

\begin{figure}[!h]
\begin{minipage}[h]{0.48\linewidth}
\center{\includegraphics[width=1\linewidth]{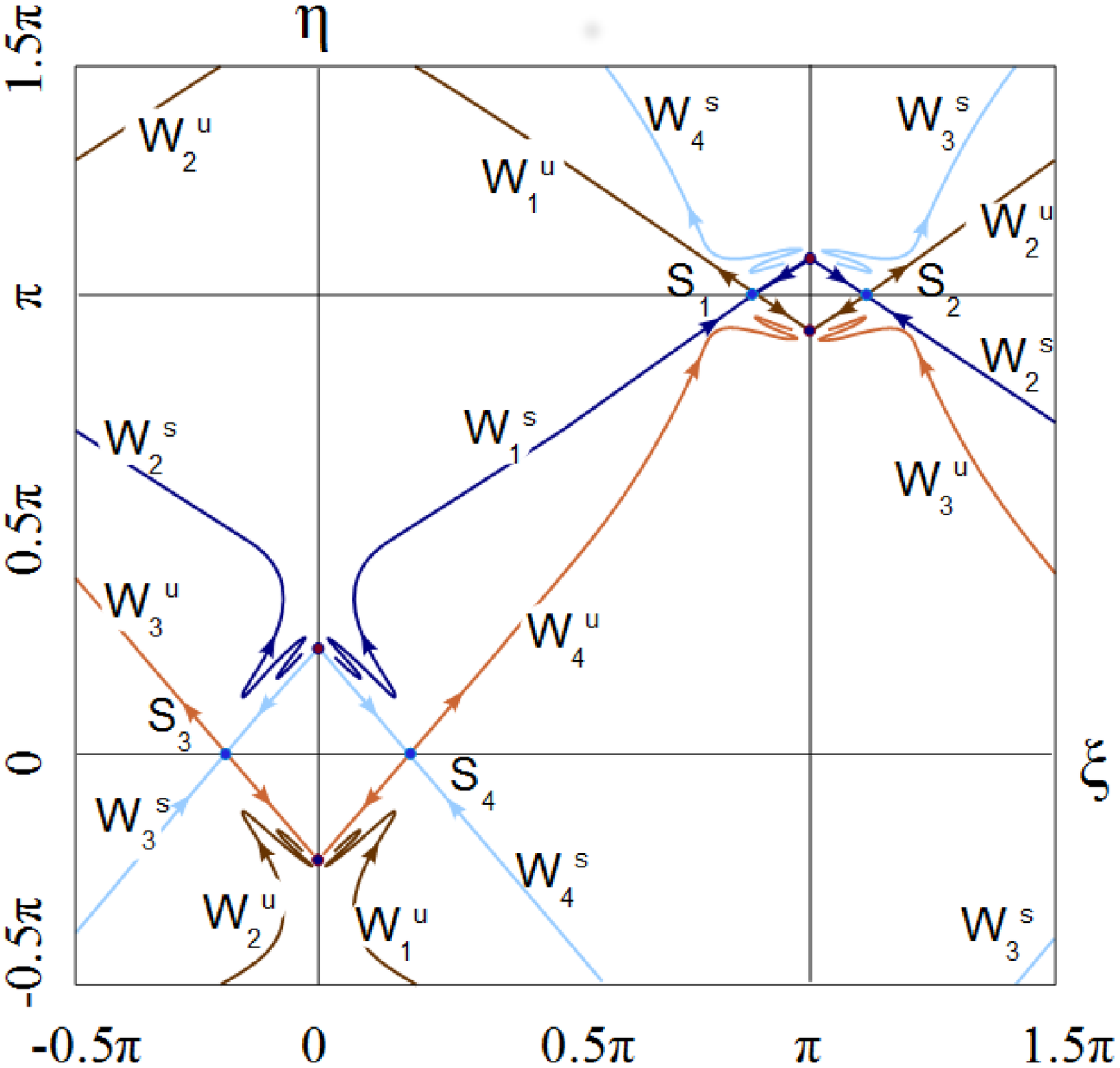} \\ (a)}
\end{minipage}
\hfill
\begin{minipage}[h]{0.48\linewidth}
\center{\includegraphics[width=1\linewidth]{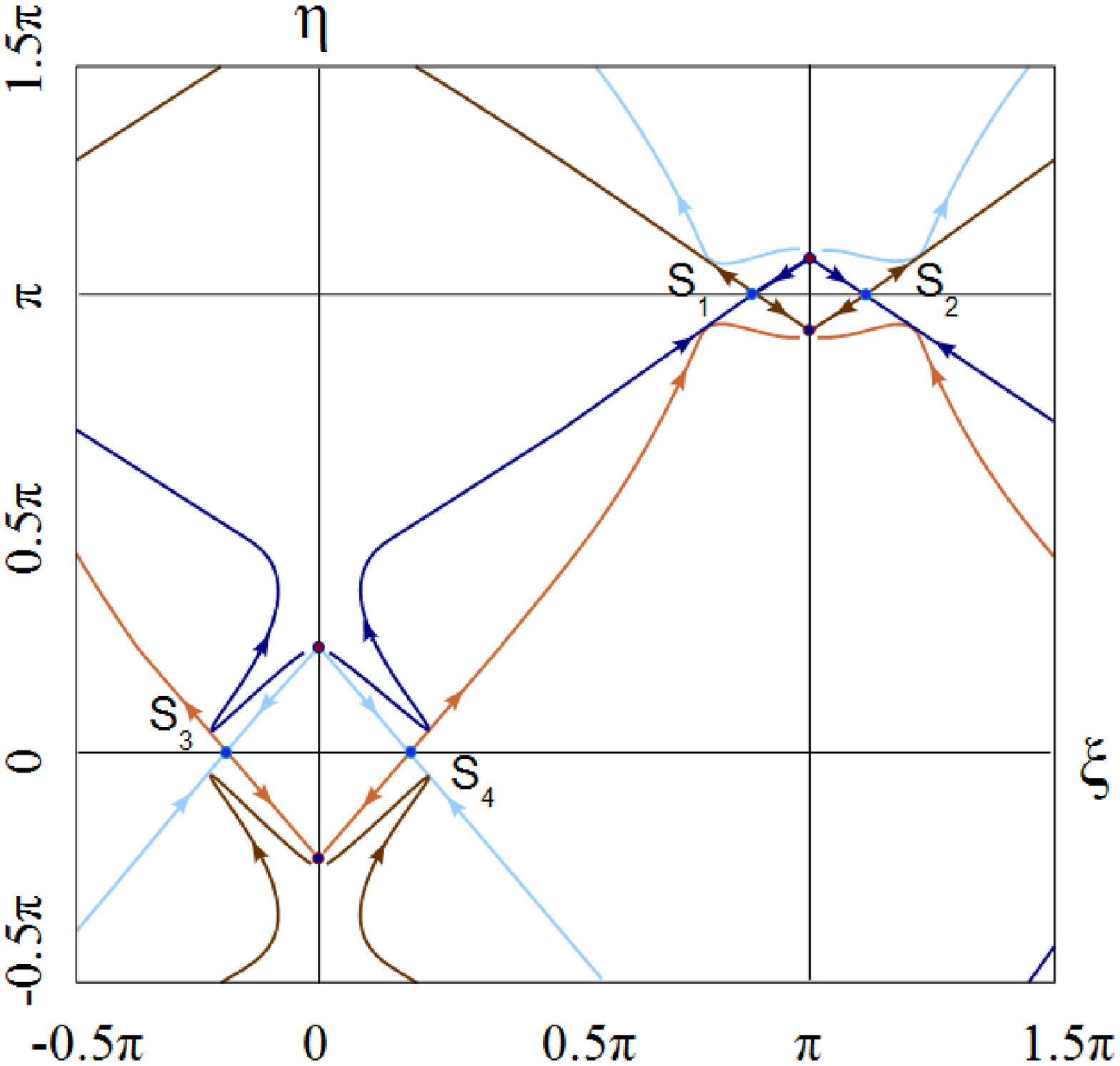} \\ (b)}
\end{minipage}
\vfill
\begin{minipage}[h]{0.48\linewidth}
\center{\includegraphics[width=1\linewidth]{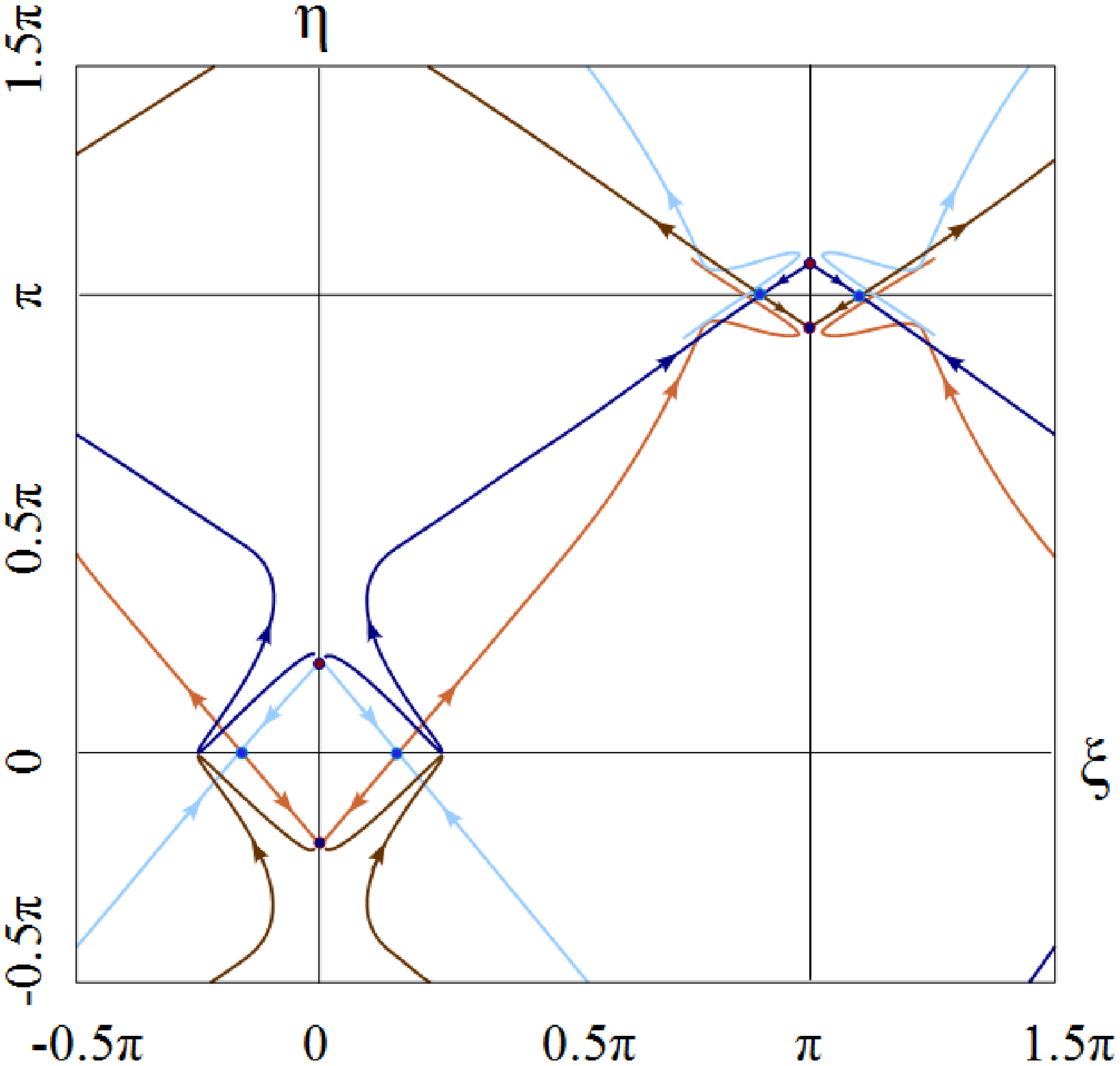} \\ (c)}
\end{minipage}
\hfill
\begin{minipage}[h]{0.48\linewidth}
\center{\includegraphics[width=1\linewidth]{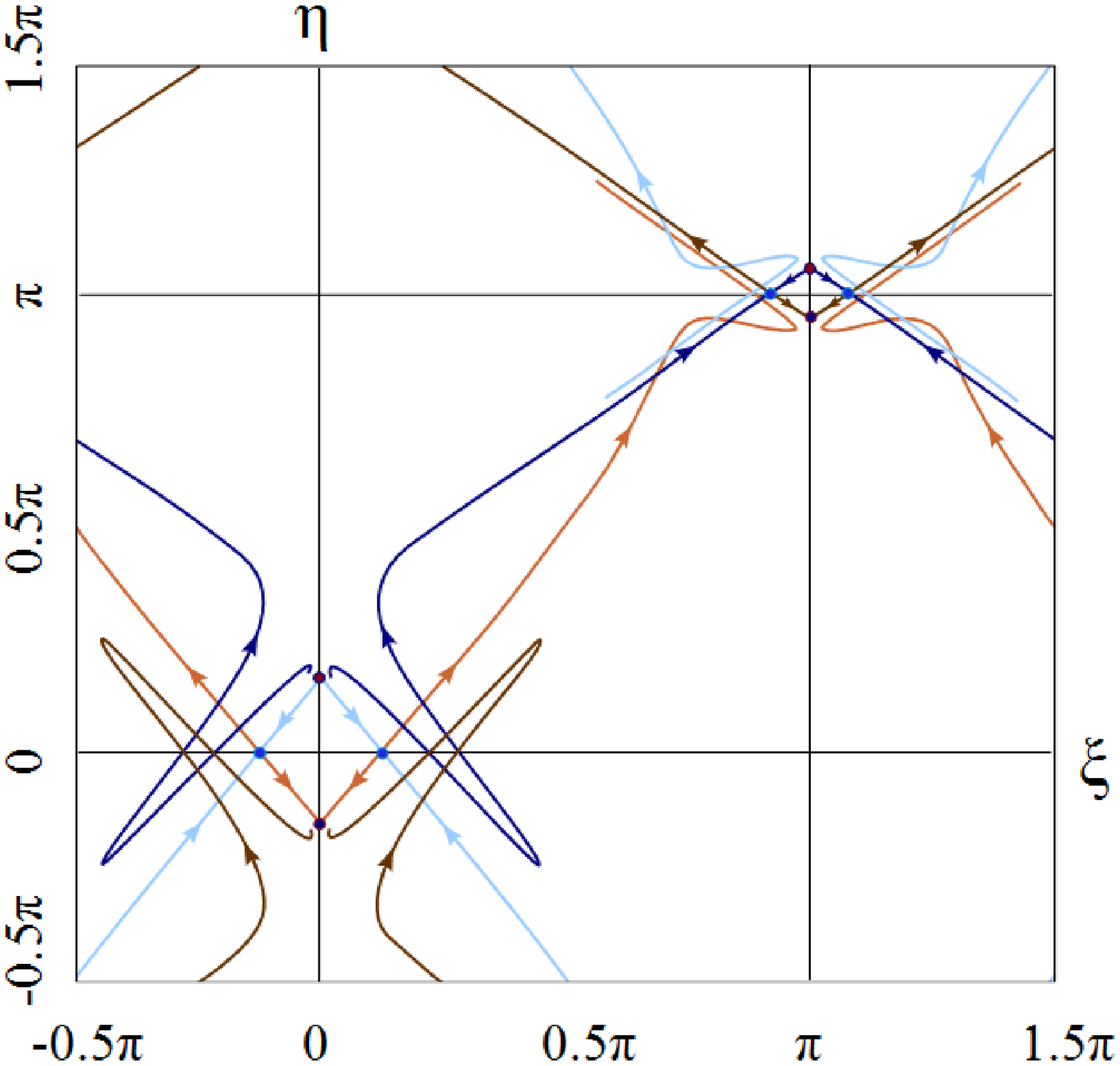} \\ (d)}
\end{minipage}
\caption{{\footnotesize The stable and unstable
manifolds of saddle fixed points of $T$ are shown. These manifolds
(a) do not intersect, $\varepsilon = 0.7$,
(b) form heteroclinic intersections, $\varepsilon = \varepsilon_1^{het} \approx 0.690$, and dynamics become chaotic.
(c) Symmetric heteroclinic orbits appear at $\varepsilon = \varepsilon_2^{het} \approx 0.679$.
(d) Developed homoclinic and heteroclinic tangles (shown at $\varepsilon = 0.650$) exist at $\varepsilon < \varepsilon_2^{het}$.}}
\label{Fig:FP1_heteroclinic_MD}
\end{figure}

The end of the mixed dynamics at larger $\varepsilon$ is related with the above described bifurcation of the birth of the fixed points
of $T$ at $\varepsilon = \varepsilon_1^*\simeq 0.6042$. Just after this bifurcation, forward iterations of $T$ for the majority of points of the phase space
seem to tend to the stable fixed points. These points become global attractors (and their images by $R$, the unstable nodes, become global repellers)
at $\varepsilon > \varepsilon_1^{het} \simeq  0.690$ when the invariant manifolds of the saddle fixed points stop having intersections.
At these values of $\varepsilon$ the dynamics of $T$ is very simple, see Fig.~\ref{Fig:FP1_heteroclinic_MD}(a). The map $T$ has 8
fixed points: four saddles $S_1,S_2,S_3$ and $S_4$, two unstable nodes (global repellers), and two stable nodes (global attractors).
In Fig.~\ref{Fig:FP1_heteroclinic_MD} we show also a picture of the stable and unstable invariant
manifolds $W^{s,u}_i$ of the saddles $S_i$, $i=1,...,4$. At $\varepsilon > \varepsilon_1^{het}$ these manifolds do not intersect and form boundaries
of the domains of attraction ($W^s_i$) and repulsion ($W^u_i$) for the global attractors and repellers.

\begin{figure}[!h]
\begin{minipage}[h]{0.48\linewidth}
\center{\includegraphics[width=1\linewidth]{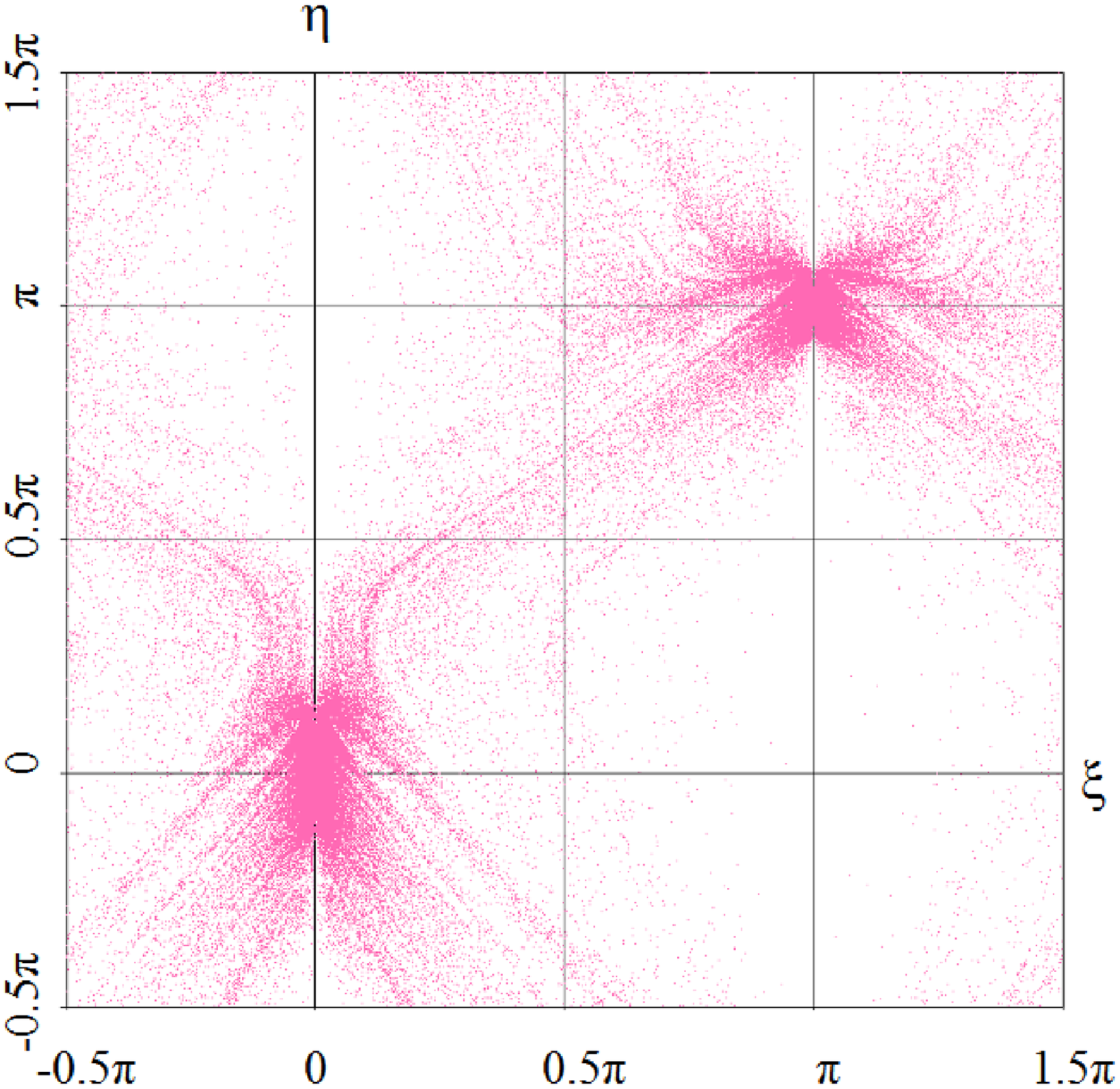} \\ (a)}
\end{minipage}
\hfill
\begin{minipage}[h]{0.48\linewidth}
\center{\includegraphics[width=1\linewidth]{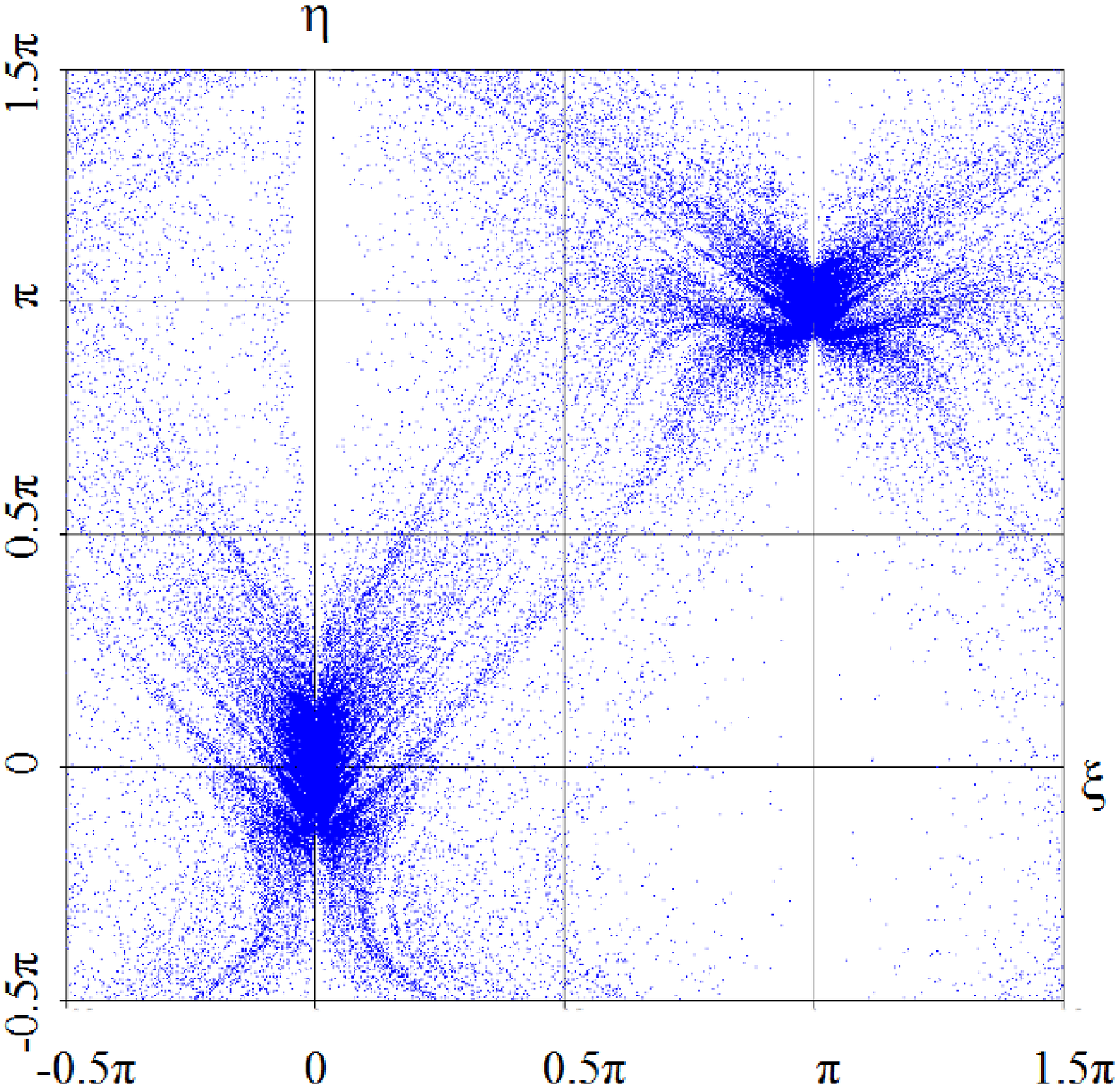} \\ (b)}
\end{minipage}
\caption{{\footnotesize  (a) Attractor and (b) repeller for $\varepsilon = 0.6 < \varepsilon_1^* = 0.6042$. One can see that the
attractor and repeller intersect but do not coincide.}}
\label{Fig:FP1_heteroclinic_MD2}
\end{figure}

At $\varepsilon < \varepsilon_1^{het} \simeq  0.690$ the manifolds $W^s_1$ and $W^u_4$, $W^s_2$ and $W^u_3$, $W^s_3$ and $W^u_2$, and $W^s_4$ and $W^u_1$ begin intersect
and heteroclinic cycles form (see Fig.~\ref{Fig:FP1_heteroclinic_MD}(b)). As a result, the dynamics becomes chaotic.
Moreover, one should associate this chaos with the mixed dynamics, since the bifurcations of these cycles such cycles lead to the birth of attractors and repellers
\cite{DGGLS13}. More heteroclinic connections emerge with the decrease of $\varepsilon$, see Figs.~\ref{Fig:FP1_heteroclinic_MD}c,d, so the chaotic set grow.
Theoretically, this chaotic set does not give a purely transient regime, since it must contain attractors and repellers born from homoclinic tangencies.
However, the numerically observed orbits do not seem to see these attractors, and converge to the stable fixed points. This coexistence of the mixed-dynamics type chaos
and stable fixed points continues until the moment $\varepsilon = \varepsilon_1^*$ when the fixed points disappear - then the large chaotic set becomes a visible
limit set for the numerically obtained trajectories. This chaos can definitely be identified with the mixed dynamics, since the attractor and repeller
intersect but do not coincide; see Fig.~\ref{Fig:FP1_heteroclinic_MD2}.

\section{Conclusion}

One can distinguish three forms of dynamical chaos: strange attractors, conservative chaos, and mixed dynamics. The
first two types are well known. Conservative chaos was, in point of fact, discovered by Poincare; the theory of strange attractors
stems from the famous work by Lorenz \cite{Lorenz}. The idea of mixed dynamics as a new type of chaos emerged
quite recently \cite{GST97}. Its distinctive feature is an unbreakable intersection of the attractor and repeller. A proper
mathematical concept of the mixed dynamics can be built with the help of the $\varepsilon$-orbit construction, going back
to works of Anosov, Conley, and Ruelle (see Section 2).

In our opinion, the mixed dynamics phenomenon must be typical for time-reversible systems \cite{DGGLS13, LSt04}
Thus, mixed dynamics is observed in systems from application. In particular, many models of nonholonomic mecanics
demonstrate such type of orbit behavior: non-holonomic constraints, though dissipation free, typically destroy the Hamiltonian structure,
but the reversibility typically persists. Thus, non-holonomic constraints produce reversible systems without a smooth invariant measure,
which is a natural condition for the emergence of mixed dynamics. For example, reversible mixed dynamics was observed in a Celtic stone model \cite{GGK_RCD_2013}, in the model of rubber Chaplygin top \cite{Kaz_RCD_2013, Kaz_NDMR_2014}, etc. 
Another, non-mechanical application is given by the Pikovsky-Topaj system of coupled rotators, which exhibits many important basic features
of the reversible mixed dynamics. With the example of this system, we show how the mixed dynamics emerges
within what appears to be a conservative chaos.
We demonstrate that the main role here is played by local and global symmetry-breaking bifurcations. Thus, local bifurcations lead to the emergence of periodic orbits of dissipative
types within the ``chaotic sea''. These orbits can be a pair of a sink and source, so the conservativity is violated
explicitly. However, there is a more subtle mechanism of the conservativity break-down, when the local symmetry-breaking bifurcation leads to the birth of
a pair of saddle periodic orbits, for which the Jacobian of the period map is less than $1$ at one of the saddles and greater than $1$ at the other one.
Invariant manifolds of these saddles intersect; moreover, as parameters change, tangencies between these manifolds also appear inevitably. In this way,
non-transverse heteroclinic cycles emerge, bifurcations of which are known to lead to the birth of infinitely many periodic attractors, repellers, and elliptic orbits;
these orbits coexist and the closures of the set of orbits of each type have non-empty intersections with each other \cite{DGGLS13,LSt04}.

Another interesting phenomenon we have found in this model is the possibility of almost instantaneous transition from a simple attractor-
repeller pair to a fully developed mixed dynamics. This corresponds to a bifurcation where the attracting fixed point, repelling fixed point, and a pair of symmetric saddles
collide (to form a degenerate saddle) and disappear. At the moment of such bifurcation, the stable and unstable separatrices of the degenerate saddle fixed
point can have homoclinic intersections, so this point can be a part of a ``large'' hyperbolic set.
The significant portion of this set survives the bifurcation, its stable and unstable manifolds can naturally acquire homoclinic tangencies as parameters change,
and this leads to mixed dynamics. Indeed, as we have seen in the Pikovsky-Topaj model, the attractor and repeller which emerge after this bifurcation
intersect over a set covering a large portion of the phase space.

\section{Acknowledgement}

This paper was supported by grant 14-41-00044 of the RSF. Numerical experiments in Section 5 were supported by
grant RSF 14-12-00811. The authors also acknowledge the support from the Royal Society grant IE141468 and
grants RFBR 14-01-00344 and 16-01-00324. AG was partially
supported by the Program of Basic Financing by the Russian Ministry of Education and Science,
AK was supported by the Basic Research Program at the National Research University
Higher School of Economics (project 98) and by the Dynasty Foundation.

\end{document}